\documentstyle[twoside]{article}

\input{amssymb.sty}
\input{times.sty}

\catcode`\@=11
\long\def\@makefntext#1{
\protect\noindent \hbox to 3.2pt {\hskip-.9pt
$^{{\eightrm\@thefnmark}}$\hfil}#1\hfill}

\def\@makefnmark{\hbox to 0pt{$^{\@thefnmark}$\hss}}	
\def\ps@myheadings{\let\@mkboth\@gobbletwo \def\@oddhead{\hbox{}
\rightmark\hfil\eightrm\thepage}
\def\@oddfoot{}\def\@evenhead{\eightrm\thepage\hfil
\leftmark\hbox{}}\def\@evenfoot{}
\def\sectionmark##1{}\def\subsectionmark##1{}}

\oddsidemargin=\evensidemargin
\newcounter{sectionc}
\newcounter{subsectionc}
\newcounter{subsubsectionc}
\renewcommand{\section}[1] {\vspace{12pt}\addtocounter{sectionc}{1}
\setcounter{subsectionc}{0}\setcounter{subsubsectionc}{0}\noindent
{\tenbf\thesectionc. #1}\par\vspace{5pt}}
\renewcommand{\subsection}[1] {\vspace{12pt}\addtocounter{subsectionc}{1}
\setcounter{subsubsectionc}{0}\noindent
{\bf\thesectionc.\thesubsectionc. {\kern1pt \bfit #1}}\par\vspace{5pt}}
\renewcommand{\subsubsection}[1]
{\vspace{12pt}\addtocounter{subsubsectionc}{1}
\noindent{\tenrm\thesectionc.\thesubsectionc.\thesubsubsectionc. {\kern1pt
\tenit #1}}\par\vspace{5pt}}
\newcommand{\nonumsection}[1] {\vspace{12pt}\noindent{\tenbf #1}
\par\vspace{5pt}}

\newcounter{appendixc}
\newcounter{subappendixc}[appendixc]
\newcounter{subsubappendixc}[subappendixc]
\renewcommand{\thesubappendixc}{\Alph{appendixc}.\arabic{subappendixc}}
\renewcommand{\thesubsubappendixc}
{\Alph{appendixc}.\arabic{subappendixc}.\arabic{subsubappendixc}}

\renewcommand{\appendix}[1] {\vspace{12pt}
\refstepcounter{appendixc}
\setcounter{figure}{0}
\setcounter{table}{0}
\setcounter{lemma}{0}
\setcounter{theorem}{0}
\setcounter{corollary}{0}
\setcounter{definition}{0}
\setcounter{equation}{0}
\renewcommand{\thefigure}{\Alph{appendixc}.\arabic{figure}}
\renewcommand{\thetable}{\Alph{appendixc}.\arabic{table}}
\renewcommand{\theappendixc}{\Alph{appendixc}}
\renewcommand{\thelemma}{\Alph{appendixc}.\arabic{lemma}}
\renewcommand{\thetheorem}{\Alph{appendixc}.\arabic{theorem}}
\renewcommand{\thedefinition}{\Alph{appendixc}.\arabic{definition}}
\renewcommand{\thecorollary}{\Alph{appendixc}.\arabic{corollary}}
\renewcommand{\theequation}{\Alph{appendixc}.\arabic{equation}}
\noindent{\tenbf Appendix#1}\par\vspace{5pt}} \newcommand{\subappendix}[1]
{\vspace{12pt}
\refstepcounter{subappendixc}
\noindent{\bf Appendix \thesubappendixc. {\kern1pt \bfit #1}} \par\vspace{5pt}}
\newcommand{\subsubappendix}[1] {\vspace{12pt}
\refstepcounter{subsubappendixc}
\noindent{\rm Appendix \thesubsubappendixc. {\kern1pt \tenit #1}}
\par\vspace{5pt}}

\newcommand{\textlineskip}{\baselineskip=13pt}
\newcommand{\smalllineskip}{\baselineskip=10pt}

\newcommand{\copyrightheading}[1]
{\vspace*{-2.5cm}\smalllineskip{\flushleft {\footnotesize Communications in
Contemporary Mathematics \\
#1}\\
{\footnotesize \copyright\kern2pt World Scientific Publishing
Company}\\
}}

\def\abstracts#1#2#3{{
\centering{\begin{minipage}{4.5in}\footnotesize\baselineskip=10pt
\parindent=0pt #1\par
\parindent=15pt #2\par
\parindent=15pt #3
\end{minipage}}\par}}

\renewenvironment{thebibliography}[1]
{\frenchspacing
\ninerm\baselineskip=11pt
\begin{list}{\arabic{enumi}.}
{\usecounter{enumi}\setlength{\parsep}{0pt}
\setlength{\leftmargin 12.7pt}{\rightmargin 0pt} 
	\setlength{\leftmargin 17pt}{\rightmargin 0pt} 
\setlength{\itemsep}{0pt} \settowidth
{\labelwidth}{#1.}\sloppy}}{\end{list}}

\newcounter{itemlistc}
\newcounter{romanlistc}
\newcounter{alphlistc}
\newcounter{arabiclistc}

\newcommand{\fcaption}[1]{
\refstepcounter{figure}
\setbox\@tempboxa = \hbox{\footnotesize Fig.~\thefigure. #1} \ifdim
\wd\@tempboxa > 5in
{\begin{center}
\parbox{5in}{\footnotesize\smalllineskip Fig.~\thefigure. #1}
\end{center}}
\else
{\begin{center}
{\footnotesize Fig.~\thefigure. #1}
\end{center}}
\fi}

\newcommand{\tcaption}[1]{
\refstepcounter{table}
\setbox\@tempboxa = \hbox{\footnotesize Table~\thetable. #1} \ifdim
\wd\@tempboxa > 5in
{\begin{center}
\parbox{5in}{\footnotesize\smalllineskip Table~\thetable. #1}
\end{center}}
\else
{\begin{center}
{\footnotesize Table~\thetable. #1}
\end{center}}
\fi}

\def\@citex[#1]#2{\if@filesw\immediate\write\@auxout
{\string\citation{#2}}\fi
\def\@citea{}\@cite{\@for\@citeb:=#2\do
{\@citea\def\@citea{,}\@ifundefined
{b@\@citeb}{{\bf ?}\@warning
{Citation `\@citeb' on page \thepage \space undefined}} {\csname
b@\@citeb\endcsname}}}{#1}}

\newif\if@cghi
\def\cite{\@cghitrue\@ifnextchar [{\@tempswatrue
\@citex}{\@tempswafalse\@citex[]}}
\def\citelow{\@cghifalse\@ifnextchar [{\@tempswatrue
\@citex}{\@tempswafalse\@citex[]}}
\def\@cite#1#2{{$\null^{#1}$\if@tempswa\typeout
{IJCGA warning: optional citation argument ignored: `#2'} \fi}}

\def\pmb#1{\setbox0=\hbox{#1}
\kern-.025em\copy0\kern-\wd0
\kern.05em\copy0\kern-\wd0
\kern-.025em\raise.0433em\box0}

\def\fnt#1#2{\footnotetext{\kern-.3em
{$^{\mbox{\scriptsize #1}}$}{#2}}}

\def\fpage#1{\begingroup
\voffset=.3in
\thispagestyle{empty}\begin{table}[b]\centerline{\footnotesize #1}
\end{table}\endgroup}

\def\runninghead#1#2{\pagestyle{myheadings}
\markboth{{\protect\footnotesize\it{\quad #1}}\hfill}
{\hfill{\protect\footnotesize\it{#2\quad}}}} \headsep=15pt

\font\tenrm=cmr10
\font\tenit=cmti10
\font\tenbf=cmbx10
\font\bfit=cmbxti10 at 10pt
\font\ninerm=cmr9

\font\eightrm=cmr8

\newtheorem{thm}{Theorem}[sectionc]

\newtheorem{cor}[thm]{Corollary}
\newtheorem{conj}{Conjecture}
\newtheorem{lemma}[thm]{Lemma}
\newtheorem{prop}[thm]{Proposition}

\newtheorem{rems}[thm]{Remarks}
\newtheorem{examples}[thm]{Examples}

\newcommand{\C}{{\Bbb C}}
\newcommand{\E}{\widetilde{\cal E}}
\newcommand{\F}{\widetilde{\cal F}}
\newcommand{\M}{\widetilde{M}}
\newcommand{\W}{\widetilde{W}}
\newcommand{\J}{{\bf H}}

\newcommand{\npartial}{\not\!\partial}
\newcommand{\SO}{{\mathbf SO}}
\newcommand{\PSL}{{\mathbf PSL}}

\newcommand{\Index}{{\em index}}

\newcommand{\rank}{{\em rank}}
\newcommand{\tr}{{\em tr}}

\newcommand{\R}{{\Bbb R}}
\newcommand{\Z}{{\Bbb Z}}
\newcommand{\Q}{{\Bbb Q}}

\newcommand{\nc}{\newcommand}
\nc{\nt}{\newtheorem}
\nc{\gf}[2]{\genfrac{}{}{0pt}{}{#1}{#2}} \nc{\mb}[1]{{\mbox{$ #1 $}}}
\nc{\real}{{\R}}
\nc{\comp}{{\C}}
\nc{\ints}{{\bf Z}}
\nc{\Ltoo}{\mb{L^2({\mathbf H})}}
\nc{\rtoo}{\mb{{\mathbf R}^2}}
\nc{\slr}{{\mathbf {SL}}(2,\real)}
\nc{\slz}{{\mathbf {SL}}(2,\ints)}
\nc{\su}{{\mathbf {SU}}(1,1)}
\nc{\SU}{{\mathbf {SU}}}
\nc{\so}{{\mathbf {SO}}}
\nc{\hyp}{{\bf H}}
\nc{\disc}{{\mathbf D}}
\nc{\torus}{{\bf T}}

\nc{\ca}{{\cal A}}
\nc{\cag}{{{\cal A}^\Gamma}}
\nc{\cg}{{\cal G}}
\nc{\chh}{{\cal H}}
\nc{\ck}{{\cal B}}
\nc{\cl}{{\cal L}}
\nc{\cm}{{\cal M}}
\nc{\cs}{{\cal S}}
\nc{\cz}{{\cal Z}}

\nc{\sind}{\sigma{\rm -ind}}

\def\qed{\hbox{${\vcenter{\vbox{	
\hrule height 0.4pt\hbox{\vrule width 0.4pt height 6pt \kern5pt\vrule width
0.4pt}\hrule height 0.4pt}}}$}}

\def\theequation{\thesectionc.\arabic{equation}}	

\begin{document}
\setlength{\textheight}{7.7truein} 

\runninghead{Twisted Index Theory on Good Orbifolds} {M. Marcolli and V.
Mathai}

\normalsize\textlineskip
\thispagestyle{empty}
\setcounter{page}{1}

\vspace*{0.88truein}

\fpage{1}
\centerline{\bf TWISTED INDEX THEORY ON GOOD ORBIFOLDS, I:} \baselineskip=13pt
\centerline{\bf NONCOMMUTATIVE BLOCH THEORY} \vspace*{0.37truein}
\centerline{\footnotesize MATILDE MARCOLLI} \baselineskip=12pt
\centerline{\footnotesize\it Department of Mathematics, Massachusetts
Institute of Technology, Cambridge, MA 02139 USA} \baselineskip=10pt
\centerline{\footnotesize\it email: matilde@math.mit.edu} \vspace*{10pt}
\centerline{\footnotesize VARGHESE MATHAI} \baselineskip=12pt
\centerline{\footnotesize\it Department of Mathematics, University of
Adelaide, Adelaide 5005, Australia}
\baselineskip=10pt
\centerline{\footnotesize\it email: vmathai@maths.adelaide.edu.au}
\vspace*{0.225truein}

\vspace*{0.21truein}
\abstracts{We study the twisted index theory of elliptic operators on
orbifold covering
spaces of compact good orbifolds, which are invariant under a projective action
of the orbifold fundamental group. We apply these results to obtain
qualitative results on real and complex hyperbolic spaces in 2 and 4
dimensions,
related to generalizations of the {\em Bethe-Sommerfeld conjecture} and the
{\em Ten Martini Problem}, on the spectrum of self adjoint elliptic
operators which are invariant under a projective action of a discrete
cocompact group.}{}{}

\vspace*{1pt}\textlineskip
\nonumsection{Introduction}
\vspace*{-0.5pt}

\noindent

In this paper, we prove a twisted index theorem for elliptic operators on
orbifold covering spaces of compact good orbifolds, which are invariant
under a projective
action of the orbifold fundamental group.

Let $\Gamma$ be a Fuchsian group of signature $(g, \nu_1,\ldots, \nu_n)$
(cf. section 1 for more details), that is, $\Gamma$ is the orbifold
fundamental group of the 2 dimensional hyperbolic orbifold $$\Sigma(g,
\nu_1,\ldots,\nu_n)$$ of signature $(g, \nu_1,\ldots,\nu_n)$. Using a result
of Kasparov \cite{Kas1} on $K$-amenable groups as well as a calculation by
Farsi \cite{Far} of the {\em orbifold} $K$-theory of compact 2-dimensional
hyperbolic orbifolds, we are able to compute the $K$-theory of twisted
group $C^*$ algebras, under the assumption that the Dixmier-Douady
invariant of the multiplier $\sigma$ is trivial
$$
K_j(C^*(\Gamma, \sigma)) \cong \left\{\begin{array}{l} \Z^{2-n
+\sum_{i=1}^n \nu_j}\qquad \hbox{if}\,\, j=0;\\ \\ {\Z}^{2g} \qquad
\qquad\qquad\hbox{if}\,\, j=1. \end{array}\right. $$
Notice that $K_0$ is much larger in the general Fuchsian group case than in
the torsion free case, where $K_0$ was determined \cite{CHMM} to be always
$\Z^2$. We also show that the orbifold $K$-theory of any 2-dimensional
orbifold is generated by orbifold line bundles.
The result is derived by means of equivariant $K$-theory and the
Baum-Connes \cite{BC} equivariant Chern character with values
in the delocalized equivariant cohomology of the smooth surface
$\Sigma_{g'}$ that covers the good orbifold $\Sigma(g,
\nu_1,\ldots,\nu_n)$. We show that the Seifert invariants \cite{Sc}
correspond to the pairing of the equivariant Chern character \cite{BC} with
a fundamental class in the delocalized equivariant homology of
$\Sigma_{g'}$.

Let $\tr$ denote the canonical trace on the twisted group $C^*$-algebra,
$C^*(\Gamma, \sigma)$, which induces a map $[\tr]$ on $K$-theory. Using the
results above and our twisted index theorem for orbifolds, we compute in
section
2 under the same assumptions as before, the range of the trace on
$K$-theory to be, $$ [\tr]( K_0 (C^*(\Gamma,\sigma))) = \Z\theta + \Z +
\sum_{i=1}^n \Z (1/\nu_i) $$ where $\theta$ denotes the evaluation of the
multiplier $\sigma$ on
the fundamental class of $\Gamma$. We then apply our calculation of the
range of
the trace on $K$-theory to the study of some quantitative aspects of the
spectrum
of projectively periodic elliptic operators on the hyperbolic plane, what
is known as {\em noncommutative Bloch theory}. Some of the most outstanding
open problems
about magnetic Schr\"odinger operators or Hamiltonians on Euclidean space
are concerned with the nature of their spectrum, most notably the {\em
Bethe-Sommerfeld conjecture} (BSC) and the {\em Ten Martini Problem} (TMP)
\cite{Sh}. More precisely, TMP asks whether given a multiplier $\sigma$ on
${\Z}^2$, is there an associated Hamiltonian ({\it i.e.} a Hamiltonian which
commutes with the $({\Z}^2, \bar\sigma)$ projective action of ${\Z}^2$ on
$L^2({\R}^2)$) possessing a Cantor set type spectrum, in the sense that the
intersection of the spectrum of the Hamiltonian with some compact interval
in ${\R}$ is a Cantor set? One can deduce from the range of the trace on
$K_0$ of the twisted group $C^*$-algebras that when the multiplier takes
its values in the roots of unity in ${\mathbf U}(1)$ (we say then that it
is rational) that
such a Hamiltonian cannot exist. However, in the Euclidean case and for
almost all irrational numbers, the discrete form of TMP has been settled in
the affirmative
\cite{Last}. BSC asserts that if the multiplier is trivial, then the
spectrum of
any associated Hamiltonian has only a {\em finite} number of gaps. This was
first
established in the Euclidean case by Skrigonov \cite{Skri}. In Sections 2
and 3,
we are concerned also with generalizations of the TMP and the BSC, which we
call
the {\em Generalized Ten Dry Martini Problem} and the {\em Generalized
Bethe-Sommerfeld conjecture}. We prove that the Kadison constant of the
twisted group $C^*$-algebra $C^*_r(\Gamma,\sigma)$ is positive whenever the
multiplier is
rational, where $\Gamma$ is now the orbifold fundamental group of a
signature $(g, \nu_1, \ldots ,\nu_n)$ hyperbolic orbifold. We then use the
results of Br\"uning and Sunada \cite{BrSu} to deduce that when the
multiplier is rational,
the generalized Ten Dry Martini Problem is answered in the negative, and we
leave
open the more difficult irrational case. More precisely, we show that the
spectrum of such a $(\Gamma, \bar\sigma)$ projectively periodic elliptic
operator
is the union of countably many (possibly degenerate) closed intervals which
can only accumulate at infinity. This also gives evidence that the
generalized Bethe-Sommerfeld conjecture is true, and generalizes earlier
results \cite{CHMM} in the torsion-free case. In order to show that the
results are not a
purely two-dimensional phenomenon, we present similar results on real and
complex hyperbolic four-manifolds, see also \cite{M}. In section 3, we
again use the range of the trace theorem above,
together with other geometric arguments to give a complete classification
up to isomorphism of the twisted group $C^*$ algebras
$C^*_r(\Gamma,\sigma)$, where $\sigma$ is assumed to have trivial
Dixmier-Douady invariant as before.

In a forthcoming paper we shall generalize these results by proving a
twisted higher index theorem which adapts the index theorems of Atiyah,
Connes and Moscovici, and Gromov, to the case of good orbifolds. This will
allow us to
compute the range of some higher traces on $K$-theory. More precisely,
suppose that $c$ is the area 2-cocycle on the group $\Gamma$ as given
above, and $\tr_c$ is the
induced cyclic 2-cocycle on a smooth subalgebra of the twisted group
$C^*$-algebra
$C^*(\Gamma, \sigma)$, which induces a map $[\tr_c]$ on $K$-theory. Then we
will prove that
$$
[\tr_c]( K_0 (C^*(\Gamma,\sigma))) = \phi\Z $$ where $\phi =
2(g-1)+ (n-\nu) \in {\Q}, \quad \nu = \sum_{j=1}^n 1/\nu_j$. By relating
the hyperbolic Connes-Kubo cyclic 2-cocycle and the area cyclic 2-cocycle,
the range of the higher trace on $K$-theory can be used to compute the
values of the Hall conductance: The results will be applied to the study
the occurrence of {\em fractional quantum numbers} in the Quantum Hall
Effect on the hyperbolic plane. We will also establish the noncommutative
Bloch theory
results for discrete Harper type operators, which is the analogue in the
discrete
case of results proved in the present paper. The results contained in this
paper, together with the results on the fractional quantum numbers, were
circulated as a preprint in 1998 \cite{MM}.

\section{Preliminaries}

\subsection{Good orbifolds}

Further details on the fundamental material on orbifolds can be found in
several references \cite{Sc}, \cite{FuSt}, \cite{Bro}. The definition of an
orbifold generalizes that of a manifold. More precisely, an {\em orbifold}
$M$ of dimension $m$ is a Hausdorff, second countable topological space
with a Satake atlas ${\cal V} = \{U_i, \phi_i\}$ which covers $M$,
consisting of open sets $U_i$ and homeomorphisms $\phi_i : U_i \to
D^m/G_i$, where $D^m$ denotes the unit ball in ${\R}^m$ and $G_i$ is a
finite subgroup of the orthogonal group $O(m)$, satisfying the following
compatibility relations; the compositions $$ \phi_j\circ\phi_i^{-1} :
\phi_i(U_i\cap U_j) \to \phi_j(U_i\cap U_j) $$ locally lifts to be a smooth
map $\R^m \to \R^m$, whenever the intersection $U_i\cap U_j \ne \emptyset$.
The open
sets $U_i$ are called local orbifold charts. In general, an orbifold $M$
can be obtained as a quotient $M=X/G$ of an infinitesimally free compact
Lie group action on a smooth manifold $X$. In fact, by Satake \cite{Sat}
and Kawasaki \cite{Kaw}, $X$ can chosen to be the smooth manifold of
orthonormal frames of the orbifold tangent bundle of $M$ (cf. section 1.4)
and $G$ can be chosen to be the orthogonal group $O(m)$.

An orbifold covering of $M$ is an orbifold map $f: Y\to M$, where $Y$ is
also an orbifold, such that any point on $M$ has a neighborhood $U$ such
that $f^{-1}(U)$ is the disjoint union of open sets $U_\alpha$, with
$f\mid_{U_\alpha}: U_\alpha \to U$ a quotient map between two quotients of
$\R^k$ by finite groups $H_1 < H_2$. The generic fibers of the covering map
$f$ are isomorphic to a discrete group which acts as deck transformations.

An orbifold $M$ is {\em good} if it is
orbifold-covered by a smooth manifold; it is {\em bad} otherwise. A good
orbifold is said to be {\em orientable} if it is orbifold covered by an
oriented manifold and the deck transformations act via orientation
preserving diffeomorphisms on the orbifold cover. Equivalently, as shown in
\cite{Sat} and \cite{Kaw}, an orbifold is orientable if it has an oriented
frame bundle $X$ such that $M = X/ \SO(m)$.

We next recall briefly some basic
notions on Euclidean and hyperbolic orbifolds, which are by fiat orbifolds
whose universal orbifold covering space is Euclidean space and hyperbolic
space respectively. We are mainly interested in the case of 2 dimensions,
and we will assume that the orbifolds in this paper are orientable.

A 2-dimensional compact orbifold has singularities that are cone points or
reflector lines. Up to passing to $\Z_2$-orbifold covers, it is always
possible to reduce to the case with only isolated cone points.

Let ${\bf H}$ denote the hyperbolic
plane and $\Gamma$ a {\em Fuchsian group of signature} $\;(g, \nu_1,\ldots,
\nu_n),$
that is, $\Gamma$ is a discrete cocompact subgroup of ${\mathbf PSL}(2,
\R)$ of genus $g$ and with $n$ elliptic elements of order $ \nu_1,\ldots,
\nu_n$ respectively. Explicitly,
$$
\begin{array}{lcl}
\Gamma & = & \Big\{A_i, B_i, C_j \in {\mathbf PSL}(2, \R) \;\Big| \;
i=1,\ldots g, \quad j = 1, \ldots n, \\[+7pt]
& & \qquad\prod_{i=1}^g[A_i, B_i] C_1\ldots C_n = 1, \quad C_j^{\nu_j} =1,
\quad j = 1, \ldots n\Big\}
\end{array}
$$
Then the corresponding compact oriented hyperbolic 2-orbifold of signature
$(g, \nu_1,\ldots, \nu_n)$ is defined as the quotient space $$ \Sigma(g,
\nu_1,\ldots,\nu_n)= \Gamma\backslash {\bf H}. $$ A compact oriented
2-dimensional Euclidean orbifold is obtained in a similar manner, but with
${\bf H}$ replaced by $\R^2$, and a complete list of these can be found in
\cite{Sc}.

Then $\Sigma(g, \nu_1,\ldots,\nu_n)$ is a compact surface of
genus $g$ with $n$ elliptic points $\{p_j\}_{j=1}^n$ such that each $p_j$
has a small coordinate neighborhood $U_{p_j} \cong D^2/{\Z_{\nu_j}}$, where
$D^2$ denotes the unit disk in $\R^2$ and $\Z_{\nu_j}$ is the cyclic group
of order ${\nu_j}, \quad j = 1, \ldots n$. Observe that the complement
$\Sigma(g, \nu_1,\ldots,\nu_n) \setminus \cup_{j=1}^n U_{p_j}$ is a compact
Riemann surface of genus $g$ and with $n$ boundary components. The group
$\Gamma$ is the orbifold fundamental group of $\Sigma(g,
\nu_1,\ldots,\nu_n)$, where the generators $C_j$ can be represented by the
$n$ boundary components of the surface $\Sigma(g, \nu_1,\ldots,\nu_n)
\setminus \cup_{j=1}^n U_{p_j}$.

All Euclidean and hyperbolic 2-dimensional orbifolds
$\Sigma(g,\nu_1,\ldots,\nu_n)$ are good, being in fact orbifold covered by
a smooth surface $\Sigma_{g'}$ cf. \cite{Sc}, i.e. there is a finite group
$G$ acting on $\Sigma_{g'}$ with quotient $\Sigma(g,\nu_1,\ldots,\nu_n)$,
where
$$ g' =1+ \frac{\#(G)}{2}(2(g -1) + (n - \nu))$$ and where $\nu =
\sum_{j=1}^n 1/\nu_j$.
According to the classification of 2-dimensional orbifolds given in
\cite{Sc}, the only bad 2-orbifolds are the ``teardrop'', with underlying
surface $S^2$ and one cone point of cone angle $2\pi/p$, and the ``double
teardrop'', with underlying surface $S^2$ and two cone points with angles
$2\pi/p$ and $2\pi/q$, $p\neq q$.

In this paper we restrict our attention to good orbifolds. It should be
pointed out that the techniques used in this paper cannot be extended
directly to the case
of bad orbifolds. It is reasonable to expect that the index theory on bad
orbifolds will involve analytical techniques for more general conic type
singularities.

\subsection{Twisted $C^*$ algebras}

We begin by recalling the definitions of the (reduced) twisted group $C^*$
algebra and its relation to a twisted $C^*$ algebra of bounded operators on
the universal orbifold cover of a good orbifold, that was defined in
\cite{BrSu}.

Let $\Gamma$ be a discrete group and $\sigma$ be a {\em multiplier} on
$\Gamma$, i.e. $\sigma:\;\Gamma\times\Gamma\to
{\mathbf U}(1)$ is a ${\mathbf U}(1)$-valued 2-cocycle on the group
$\Gamma$ i.e. $\sigma$
satisfies the following identity:
$$
\begin{array}{ll}
\bullet\quad & \sigma(\gamma,1)=\sigma(1,\gamma)=1\quad\forall \gamma\in
\Gamma;\\
\bullet\quad &\sigma(\gamma_1,\gamma_2)\sigma(\gamma_1\gamma_2,\gamma_3)=
\sigma(\gamma_1,\gamma_2\gamma_3)\sigma(\gamma_2,\gamma_3) \forall
\gamma_1,\gamma_2,\gamma_3\in \Gamma.
\end{array}
$$
Consider the Hilbert space
of square summable functions on $\Gamma$,
\[
L^2(\Gamma)=\left\{f:\Gamma\to \C\;\Big|
\;\sum_{\gamma\in\Gamma}|f(\gamma)|^2<\infty \right\}\;.
\] There are natural left $\sigma$-regular and right $\sigma$-regular
representations on $L^2(\Gamma)$. The {\em left $\sigma$-regular
representation} is defined as follows:  $\quad\forall\;\gamma,\gamma'\in
\Gamma$,
$$\begin{array}{lcl}
(L_{\gamma}^{\sigma}
f)(\gamma') & = &f(\gamma^{-1}\gamma')\sigma(\gamma,\gamma^{-1} \gamma');\\
L_{\gamma}^{\sigma} L_{\gamma'}^{\sigma} & = & \sigma(\gamma,\gamma')
L_{\gamma\gamma'}^{\sigma}.
\end{array}  $$
The {\em right $\sigma$-regular representation} is defined as follows: $\quad
\forall \gamma,\gamma'\in \Gamma$,
$$\begin{array}{lcl}
(R_{\gamma}^{\sigma} f)(\gamma') & = &
f(\gamma'\gamma)\sigma(\gamma',\gamma);\\
R_{\gamma}^{\sigma} R_{\gamma'}^{\sigma} & = &\sigma(\gamma,\gamma')
R_{\gamma\gamma'}^{\sigma}.
\end{array} $$
One can use the cocycle identity to show that the left
$\sigma$-regular representation commutes with the right $\bar{\sigma}$-regular
representation, where $\bar{\sigma}$ denotes the conjugate cocycle. Also the
left $\bar{\sigma}$-regular representation commutes with the right
$\sigma$-regular representation.
Define
\[
\qquad W^*(\Gamma,\sigma)=\left\{A\in B(\ell^2(\Gamma)):
[L_{\gamma}^{\bar\sigma}, A]=0\;
\forall \gamma\in \Gamma\right\}
\] i.e. $W^*(\Gamma,\sigma)$ is the commutant of the left
$\bar{\sigma}$-regular
representation. By general theory, it is a von Neumann algebra, and it is
called
the {\em twisted group von Neumann algebra}. It can also be realized in the
following manner: the right $\sigma$-regular representation of $\Gamma$ extends
to a $*$ representation of the twisted group algebra, $\C(\Gamma,\sigma) \to
B(\ell^2(\Gamma))$.
Now the {\em weak closure} (which coincides with the strong closure)
of $\C(\Gamma,\sigma)$ also yields the twisted group von Neumann algebra
$ W^*(\Gamma,\sigma)$, by the commutant theorem of von Neumann.
The {\em norm closure} of $\C(\Gamma,\sigma)$ yields the (reduced) {\em
twisted group $C^*$ algebra} $C^*_r(\Gamma, \sigma)$. We will also briefly
consider the full or unreduced
{\em twisted group $C^*$ algebra} $C^*(\Gamma, \sigma)$, which is another $C^*$
completion of $\C(\Gamma,\sigma)$ using $*$-representations, cf. \cite{PR}
for the definition.

Let $M$ be a good, compact orbifold, and ${\cal E} \to M$ be an orbifold
vector bundle over $M$, and $\E \to \widetilde{M}$ be its lift to the
universal orbifold covering space $\Gamma\to \widetilde{M}\to M$, which is
by assumption
a simply-connected smooth manifold.
We will now briefly review how to construct a $(\Gamma, \bar\sigma)$-action
(where $\sigma$ is a multiplier on $\Gamma$ and $\bar\sigma$ denotes its
complex conjugate) on $L^2(\M)$. Let $\omega = d\eta$ be an exact $2$-form
on $\widetilde{M}$ such that $\omega$ is also $\Gamma$-invariant, although
$\eta$ is {\em not} assumed to be $\Gamma$-invariant. Define a Hermitian
connection on the trivial line bundle over $\widetilde{M}$ as
\[
\nabla = d + i\eta.
\]
Its curvature is $\nabla^2=i\omega$. Then $\nabla$ defines a
$(\Gamma,\bar{\sigma})$ action on
$L^2(\widetilde{M},\widetilde{{\cal S}^\pm\otimes E})$ as follows:

Since $\omega$ is $\Gamma$ invariant, one has $\forall\gamma\in\Gamma$ \[ 0
= \gamma^* \omega-\omega = d(\gamma^* \eta-\eta), \] so that
$\gamma^*\eta-\eta$ is a closed $1$-form on a simply-connected manifold
$\widetilde{M}\Rightarrow \gamma^*\eta-\eta=d\psi_\gamma$ for some smooth
function $\psi_\gamma$ on $\widetilde{M}$ satisfying $$ \begin{array}{ll}
\bullet\quad & \psi_\gamma(x)+\psi_{\gamma'}(\gamma
x)-\psi_{\gamma'\gamma}(x) \qquad\hbox{is independent of $x\ \forall
x\in\widetilde{M}$};\\
\bullet\quad & \psi_\gamma(x_0)=0\qquad\hbox{for some $x_0\in\widetilde{M}
\ \ \forall \gamma\in\Gamma$}. \end{array} $$
Then
$\bar{\sigma}(\gamma,\gamma')=\exp(i\psi_\gamma(\gamma' x_0))$ defines a
multiplier on $\Gamma$.
Now define the
$(\Gamma,\bar{\sigma})$ action as follows: For $u\in L^2(\widetilde{M},
\widetilde{{\cal S}^\pm\otimes E}),
\quad U_\gamma\ u=\gamma^* u, \quad
S_\gamma u = \exp(i\,\psi_\gamma)\ u$, define $T_\gamma=U_\gamma\circ
S_\gamma$.
Then it satisfies $T_{\gamma_1} T_{\gamma_2}
=\bar{\sigma}(\gamma_1,\gamma_2)\,T_{\gamma_1 \gamma_2}$, and so it defines
a $(\Gamma,\bar{\sigma})$-action. It can be shown that only multipliers
$\bar \sigma$ such that the Dixmier-Douady invariant $ \delta(\bar\sigma) =
0$ can give
rise to $(\Gamma, \bar\sigma)$-actions in this way cf. section 2.2 for a
further discussion.

Let $D: L^2(\widetilde{M}, \E) \to L^2(\widetilde{M}, \E)$ be a self
adjoint elliptic differential
operator that commutes with a $(\Gamma, \bar\sigma)$-action $T_\gamma \quad
\forall \gamma\in \Gamma$ on $L^2(\widetilde{M}, \E)$. Then by the
functional calculus, all the spectral projections of $D$, $E_\lambda =
\chi_{[0,\lambda]}(D)$ are bounded self adjoint operators on
$L^2({\widetilde M}, \E)$ that commute with $T_\gamma \quad \forall
\gamma\in \Gamma$.
Now the commutant of the $(\Gamma, \bar\sigma)$-action on $L^2({\widetilde
M}, \E)$
is a von Neumann algebra
$$
W^*( \sigma) = \left\{Q\in B(L^2({\widetilde M}, \E))| \quad T_\gamma Q =
Q T_\gamma \quad \forall \gamma\in \Gamma \right\}. $$
Since $T_\gamma Q = Q T_\gamma $, one sees that $$
e^{-i\phi_\gamma(x)}k_Q(\gamma x,\gamma y) e^{i\phi_\gamma(y)} = k_Q(x,y)
$$ $\forall x, y \in {\widetilde M}\quad \forall \gamma\in \Gamma$, where
$k_Q$ denotes the Schwartz kernel of $Q$. In particular, observe that
$\underline\tr(k_Q(x,x))$ is a $\Gamma$-invariant function on ${\widetilde
M}$, where $\underline\tr$ denotes the pointwise trace. Using this, one
sees that there is a semi-finite trace on this von Neumann algebra $$
\tr: W^*( \sigma) \to \C
$$
defined as in the untwisted case due to Atiyah \cite{At}, $$ Q \to \int_{M}
\underline\tr(k_Q(x,x))dx
$$
where $k_Q$ denotes the Schwartz kernel of $Q$. Note that this trace is
finite whenever $k_Q$ is continuous in a neighborhood of the diagonal in
${\widetilde M}\times {\widetilde M}$.

By elliptic regularity, the spectral projection $E_\lambda$ has a smooth
Schwartz kernel, so that in particular, the spectral density function,
$N_Q(\lambda) = \tr(E_\lambda) <\infty \quad \forall \lambda$, is well
defined.

If ${\cal F}$ is a relatively compact fundamental domain in $\widetilde M$
for the action of $\Gamma$ on ${\widetilde M}$, then
one sees that there is a $(\Gamma, \bar\sigma)$-isomorphism
\begin{equation}
L^2({\widetilde M}, \E) \cong  L^2(\Gamma) \otimes L^2({\cal F}, \E|_{{\cal
F}})
\end{equation}
It is given by
$f\mapsto g$ where $ g(\gamma)(x) = f(\gamma x)$, $x\in {\cal F}, \;
\gamma\in \Gamma$,
equivalently by a choice of a bounded measurable almost everywhere smooth
section of the
orbifold covering $\widetilde M \to M$. The $(\Gamma, \bar\sigma)$-action on
$L^2({\cal F}, \E|_{{\cal F}}) $ is trivial and its is the previously
defined regular
$(\Gamma, \bar\sigma)$ representation on $L^2(\Gamma)$.
Therefore
$$
W^*( \sigma) \cong
W^*(\Gamma, \sigma) \otimes B(L^2({\cal F}, \E|_{{\cal F}}) ) $$
where $B(L^2({\cal F}, \E|_{{\cal F}}) ) $ denotes the  algebra of all
bounded operators
on the Hilbert space $L^2({\cal F}, \E|_{{\cal F}}) $.
There is a natural subalgebra ${C^*} ( \sigma)$ of $W^* ( \sigma)$ which is
defined as follows. Let
$$ \begin{array}{rl}
{C_c}^\infty ( \sigma)=
& \Big\{Q\in W^* ( \sigma)
\Big| \begin{array}{l} k_Q \ \hbox{is smooth and supported in a compact }\\
\hbox{neighborhood of the diagonal}\end{array}\Big\} \end{array}
$$
Then ${C^*}( \sigma)$ is defined to be the {\em norm closure} of
${C_c}^\infty ( \sigma)$.
It can also be shown to be the norm closure of
$$ \left\{Q\in W^* (\sigma) \Big| \begin{array}{l} k_Q \ {\hbox{is smooth
and}}\ k_Q(x,y) \ {\hbox {is $L^1$ in both}} \\ {\hbox {the $x$ and $y$
variables separately}}\end{array} \right\} $$
That is, elements of ${C^*} ( \sigma)$ are elements of $W^* (\sigma)$ that
have the additional property of
some off-diagonal decay.
Via the isomorphism given in equation $(1.1)$, it can be shown that
\begin{equation}
{C^*} ( \sigma) \cong C^*_r(\Gamma, \sigma) \otimes{\cal K}
\end{equation}
where ${\cal K} = {\cal K}(L^2({\cal F}, \E|_{\cal F}))$ denotes the $C^*$
algebra of compact operators on the Hilbert space $L^2({\cal F}, \E|_{\cal
F})$ (see
\cite{BrSu} for details).

\subsection{The $C^*$ algebra of an orbifold} Let $M$ be an oriented
orbifold of dimension $m$, that is $M = P/\SO(m)$, where $P$ is the bundle
of oriented frames on the orbifold tangent bundle (cf. section 1.4). Then
the $C^*$ algebra of the orbifold $M$ is by fiat the crossed product
$C^*(M) = C(P) \rtimes \SO(m)$, where $C(P)$ denotes the $C^*$ algebra of
continuous functions on $P$. We will now study some Morita equivalent
descriptions of $C^*(M)$ that will be useful for us later. The following is
one such, and is due to \cite{Far}.

\begin{prop} Let $M$ be a good orbifold, which is orbifold covered by the
smooth manifold $X$, i.e. $M = X/G$. Then the $C^*$ algebras $C_0(X)
\rtimes G$ and $C^*(M)$ are strongly Morita equivalent. \label{prop1}
\end{prop}

In the two dimensional case, there is yet another $C^*$ algebra that is
strongly Morita equivalent to the $C^*$ algebra of the orbifold. Let
$\Gamma$ be as before. Then $\Gamma$ acts freely on $\PSL(2, \R)$, and
therefore the quotient space
$$\Gamma\backslash \PSL(2, \R) =P(g,
\nu_1,\ldots,\nu_n)$$ is a smooth compact manifold, with a right action of
$\SO(2)$ that is only infinitesimally free. The $C^*$ algebra of the
hyperbolic orbifold $\Sigma(g, \nu_1,\ldots,\nu_n)$ is by fiat the crossed
product $C^*$ algebra
$$
C^*(\Sigma(g, \nu_1,\ldots,\nu_n)) =
C(P(g, \nu_1,\ldots,\nu_n)) \rtimes \SO(2) $$ cf. \cite{Co}. If $\SO(2)$
did act freely on $P(g, \nu_1,\ldots,\nu_n)$ (which is the case when
$\nu_1=\ldots =\nu_n=1$), then it is known that $C^*(\Sigma(g,
\nu_1,\ldots,\nu_n)) $ and $C(\Sigma(g, \nu_1,\ldots,\nu_n)) $ are strongly
Morita equivalent as $C^*$ algebras.

We shall next describe a natural algebra which is Morita equivalent to the
$C^*$ algebra of the orbifold $\Sigma(g, \nu_1,\ldots,\nu_n)$. Now $\Gamma$
has a torsion free subgroup $\Gamma_{g'}$ of finite index, such that the
quotient $\Gamma_{g'}\backslash \J = \Gamma_{g'} \backslash
\PSL(2,\R)/\SO(2) = \Sigma_{g'}$ is a compact Riemann surface of genus $g' =
1+ \frac{\#(G)}{2}(2(g -1) +(n - \nu))$
where
$\nu = \sum_{j=1}^n 1/\nu_j$, cf. Theorem 2.5 \cite{Sc}, and the orbifold
Euler characteristic calculations in there. Then $G \to \Sigma_{g'} \to
\Sigma(g, \nu_1,\ldots,\nu_n)$ is a finite orbifold cover, i.e. a ramified
covering space, where $G = \Gamma_{g'}\backslash \Gamma$. \begin{prop}
The $C^*
$ algebras $C(\Sigma_{g'})\rtimes G$, $C^*(\Sigma(g, \nu_1,\ldots,\nu_n))$
and $C_0({\bf H})\rtimes \Gamma$ are strongly Morita equivalent to each
other. \label{prop2}
\end{prop}

\noindent{\bf Proof.} The strong Morita equivalence of the last two $C^*$
algebras is contained in the previous Proposition. Since strong Morita
equivalence is an equivalence relation, it suffices to prove that the first
two $C^*$ algebras are strongly Morita equivalent. Let $$\widehat P =
\Gamma_{g'} \backslash \PSL(2,\R)$$ where $\SO(2)$ acts on $\widehat P$ the
right, and therefore commutes with the left $G$ action on $\widehat P$.
Moreover, the actions of $G$ and $\SO(2)$ on $\widehat P$ are free, and
therefore one can apply a theorem of Green, \cite{Green}
which implies in particular that $C_o(G\backslash \widehat P)\rtimes \SO(2)
$ and $C_o(\widehat P/\SO(2))\rtimes G $
are strongly Morita equivalent, i.e.
$C_o(P(g, \nu_1,\ldots,\nu_n)\rtimes \SO(2) $ and $C_o(\Sigma_{g'}) \rtimes
G $ are strongly Morita equivalent, proving the proposition. $\diamond$

\subsection{Orbifold vector bundles and $K$-theory} Because of the Morita
equivalences of the last section 1.3, we can give several alternate and
equivalent descriptions of orbifold vector bundles over orbifolds. Firstly,
there is the description using transition functions cf. \cite {Sat},
\cite{Kaw}. Equivalently, one can view an orbifold vector bundle over an
$m$ dimensional orbifold $M$ as being an $\SO(m)$ equivariant vector bundle
over the bundle $P$ of oriented frames of the orbifold tangent bundle. In
the case of a good orbifold $M$, which is orbifold covered by a smooth
manifold $X$, let $G$ be the discrete group acting on $X$, $G \to X \to
M=X/G$. Then an orbifold vector bundle on $M$ is the quotient ${\cal
V}_M=G\backslash {\cal V}_X$ of a vector bundle over $X$ by the $G$ action.
Notice that an orbifold vector bundle is not a vector bundle over $M$: in
fact, the fiber at a singular point is isomorphic to a quotient of a vector
space by a finite group action.

The Grothendieck group of isomorphism classes of orbifold vector bundles on
the orbifold $M$ is called the {\em orbifold $K$-theory} of $M$ and is
denoted by $K^0_{orb}(M)$, which by a result of \cite{BC}, \cite{Far} is
canonically isomorphic to $K_0(C^*(M))$. By the Morita equivalence of
section 1.3, one then has $K^0_{orb}(M) \cong K^0_{\SO(m)}(P)$, and by the
Julg-Green theorem \cite{Ju}, \cite{Green}, the second group is isomorphic
to $K_0(C(P) \rtimes \SO(m))$. In the case when $M$ is a good orbifold, by
Proposition \ref{prop1}, one sees that $K^0_{orb}(M) \cong K^0(C_0(X)
\rtimes G) = K^0_G (X) $.

We will now be mainly interested in orbifold line bundles over
hyperbolic 2-orbifolds. Let $G$ be the finite group determined by the exact
sequence
$$1\to \Gamma_{g'}\to \Gamma\to G\to 1.$$ Then $G$ acts on $\Sigma_{g'}$ with
quotient the orbifold $\Sigma(g,
\nu_1,\ldots,\nu_n)$.

An orbifold line bundle ${\cal L}$ on $\Sigma(g, \nu_1,\ldots,\nu_n)$ is
given by
$$ {\cal L} =G\backslash (P\times_{\SO(2)} \C),$$ where $P$ is a principal
$\SO(2)$-bundle on the smooth surface $\Sigma_{g'}$. Notice that the
$\SO(2)$ and the $G$ actions commute, and are free on the total space $P$.
An orbifold line bundle has an associated Seifert fibered space $G
\backslash P$. A more explicit local geometric construction of ${\cal L}$
is given in \cite{Sc}.
An orbifold line bundle ${\cal L}$ over a hyperbolic orbifold
$\Sigma(g,\nu_1,\ldots,\nu_n)$ is specified by the Chern class of the
pullback line bundle on the smooth surface $\Sigma_{g'}$, together with the
Seifert data. That is the pairs of numbers $(\beta_j,\nu_j)$, where
$\beta_j$ satisfies the following condition. Given the exact sequence $$
1\to \Z \to \pi_1(P)\to \pi_1^{orb}(\Sigma(g,\nu_1,\ldots,\nu_n))\to 1,$$
let $\widetilde C_j$ be an element of $\pi_1(P)$ that maps to the generator
$C_j$ of the orbifold fundamental group. Let $C$ be the generator of the
fundamental group of the fiber. Then we have $C_j^{\nu_j}=1$ and
$C^{\beta_j}=\widetilde C_j^{\nu_j}$. The choice of $\beta_j$ can be
normalized so that $0<\beta_j <\nu_j$.

More geometrically, let $\Sigma(g,\nu_1,\ldots,\nu_n)$ be a hyperbolic
orbifold with the cone points $p_1, \ldots, p_n$. Let $\Sigma'$ be the
complement of the union of small disks around the cone points. The orbifold
line bundle induces a line bundle ${\cal L}'$ over the smooth surface with
boundary $\Sigma'$, trivialized over the boundary components of $\Sigma'$.
Moreover, the restriction of the orbifold line bundle ${\cal L}$ over the
small disks $D_{p_i}$ around each cone point $p_i$ is obtained by
considering a surgery on the trivial product $\C\times D_{p_i}$ obtained by
cutting open along a radius in $\C$ and gluing back after performing a
rotation on $D_{p_i}$ by an angle $2\pi q/\nu_i$. With this notation the
Seifert invariants are $(q_i,\nu_i)$ with $\beta_i q_i\equiv 1$ $(mod
\nu_i)$.

Thus, an orbifold line bundle has a finite set of singular fibers at the
cone points. The orbifold line bundle ${\cal L}$ pulls back to a
$G$-equivariant line bundle $\widetilde {\cal L}$ over the smooth surface
$\Sigma_{g'}$ that orbifold covers $\Sigma(g,\nu_1,\ldots,\nu_n)$. All the
orbifold line bundles with trivial orbifold Euler class, as defined in
\cite{Sc}, lift to the trivial line bundle on $\Sigma_{g'}$.

In \cite{Sc} the classification of Seifert-fibered spaces is derived using
the Seifert invariants, namely the Chern class of the line bundle
$\widetilde {\cal L}$, together with the Seifert data $(\beta_j, \nu_j)$ of
the singular fibers at the cone points $p_j$. We show in the following that
the Seifert invariants can be recovered from the image of the Baum-Connes
equivariant Chern character \cite{BC}.

\subsection{Baum-Connes Chern character}

We have seen that the algebra $C^*(\Sigma(g,\nu_1,\ldots,\nu_n))$ is
strongly Morita equivalent to the cross product $C(\Sigma_{g'})\rtimes G$.
Therefore the relevant K-theory is $$K_0(C(\Sigma_{g'})\rtimes
G)=K^0_{\SO(2)}(G\backslash \widehat P)=K^0_G(\Sigma_{g'}),$$ where
$\widehat P=\Gamma_{g'} \backslash \PSL(2,\R)$.

We recall briefly the definition of delocalized equivariant cohomology for
a finite group action on a smooth manifold \cite{BC}. Let $G$ be a finite
group acting smoothly and properly on a compact smooth manifold $X$. Let
$M$ be the good orbifold $M=G\backslash X$. Given any $\gamma\in G$, the
subset
$X^\gamma$ of $X$ given by
$$ X^\gamma=\{ (x,\gamma)\in X\times G \mid \gamma x=x \} $$ is a smooth
compact submanifold. Let $\widehat X$ be the disjoint union of the
$X^\gamma$ for $\gamma\in G$. The complex $\Omega_G(\widehat X)$ of
$G$-invariant de Rham forms on $\widehat X$ with coefficients in $\C$
computes the delocalized equivariant cohomology $H^\bullet (X,G)$, which is
$\Z_2$ graded by forms of even and odd degree. The dual complex that
computes delocalized homology is
obtained by considering $G$-invariant de Rham currents on $\widehat X$.
Thus we have
$$ H^\bullet (X,G)=H^\bullet (\Omega_G(\widehat X),d)=H^\bullet (\widehat X
/G,\C)$$ $$=H^\bullet (\widehat X,\C)^G =\bigoplus_{\gamma\in G}
H^\bullet(X^\gamma,\C).$$ According to \cite{BC}, Theorem 7.14, the
delocalized equivariant cohomology is isomorphic to the cyclic cohomology
of the algebra $C^\infty(X)\rtimes G$,
$$ H^0(X, G)\cong HC^{ev}(C^\infty(X)\rtimes G), $$ $$ H^1(X, G)\cong
HC^{odd}(C^\infty(X)\rtimes G). $$

The Baum-Connes equivariant Chern character $$ ch_G : K^0_G(X)\to H^0(X,G)
$$ is an isomorphism over the complex numbers. Equivalently, the
Baum-Connes equivariant Chern character can be viewed as $$ ch_G :
K^0_{orb}(M)\to H^0_{orb}(M) $$ where the {\em orbifold cohomology} is by
definition $ H^j_{orb}(M) = H^j(X, G)$
for $j=0,1$.

In our case the delocalized equivariant cohomology and the Baum-Connes
Chern character have a simple expression. In fact, let $\Sigma_{g'}$
be the smooth
surface that orbifold covers $\Sigma(g,\nu_1,\ldots,\nu_n)$. Let $G$ be the
finite group $1\to \Gamma_{g'} \to \Gamma\to G\to 1$. Let $G_{p_j}\cong
\Z_{\nu_j}$ be
the stabilizer
of the cone point $p_j$ in $\Sigma(g,\nu_1,\ldots,\nu_n)$. Then we have $$
\Sigma_{g'}^\gamma=\left\{\begin{array}{lr} \Sigma_{g'} & \quad\hbox{if}
\quad\gamma=1; \\ \{p_j\} & \quad\hbox{if} \quad\gamma\in G_{p_j}\backslash
\{ 1 \}; \\ \emptyset & \hbox{otherwise.} \end{array}\right. $$ Thus the
delocalized equivariant cohomology and orbifold cohomology is given by $$
H^0_{orb}(\Sigma(g,\nu_1,\ldots,\nu_n)) = H^0(\Sigma_{g'},
G)=H^0(\Sigma_{g'})\oplus H^2(\Sigma_{g'})\oplus \C^{\sum_j (\nu_j-1)}, $$
where each $\C^{\nu_j-1}$ is given by $\nu_j-1$ copies of $H^0(p_j)$, and
$$ H^1_{orb}(\Sigma(g,\nu_1,\ldots,\nu_n)= H^1(\Sigma_{g'},
G)=H^1(\Sigma_{g'}). $$

Let ${\cal L}$ be an orbifold line bundle in $K_0(C(\Sigma_g)\rtimes
G)=K^0_G(\Sigma_g)$, and let $\widetilde {\cal L}$ be the corresponding
line bundle over the surface $\Sigma_{g'}$. An element $\gamma$ in the
stabilizer $G_{p_j}$ acts on the restriction of ${\cal
L}|_{\Sigma^\gamma_g}={\cal L}|_{p_j}=\C$ as multiplication by
$\lambda(\gamma)=e^{2\pi i \beta_j/\nu_j}$.

Thus, the Baum-Connes Chern character of ${\cal L}$ is given by $$
ch_G({\cal L})=(1,c_1(\widetilde {\cal L}), e^{2\pi i \beta_1/\nu_1},
\ldots, e^{2\pi i (\nu_1-1) \beta_1/\nu_1}, \ldots e^{2\pi i
\beta_n/\nu_n}, \ldots, e^{2\pi i (\nu_n-1) \beta_n/\nu_n}). $$

\begin{prop}
The Baum-Connes Chern character classifies orbifold line bundles over the
orbifold $\Sigma(g,\nu_1,\ldots,\nu_n)$. \label{prop3} \end{prop}

\noindent{\bf Proof.}
According to \cite{Sc} the orbifold line bundles are classified by the
orbifold Euler number
$$ e(\Sigma(g,\nu_1,\ldots,\nu_n))=<c_1(\widetilde {\cal L}),
[\Sigma_{g'}]>+ \sum_{j} \beta_j/\nu_j, $$ given in terms of the Chern
number $<c_1(\widetilde {\cal L}), [\Sigma_{g'}]>$ and the Seifert
invariants $(\beta_j,\nu_j)$. $\diamond$

Notice that we have the isomorphism in $K$-theory,
$K^0_G(\Sigma_{g'})=K^0_{\SO(2)}(G\backslash \widehat P)$ and the Chern
character isomorphisms (with $\C$ coefficients) $$ ch_G: K^0_G(\Sigma_{g'})
\to H^0(\Sigma_{g'}, G) \cong HC^{ev}(C^\infty(\Sigma_{g'})\rtimes G) $$
and $$ ch_{\SO(2)} : K^0_{\SO(2)}(\Gamma \backslash \PSL(2,\R))\to
HC^{ev}(C^\infty(\Gamma \backslash \PSL(2,\R))\rtimes \SO(2)).$$ Moreover,
we have an isomorphism
$$HC^\bullet (C^\infty(\Gamma \backslash \PSL(2,\R))\rtimes \SO(2))\cong
H^\bullet_{\SO(2)} (\Gamma \backslash \PSL(2,\R)). $$ Thus, we obtain
$$HC^{ev}(C^\infty(\Gamma \backslash \PSL(2,\R))\rtimes \SO(2)) \cong
HC^{ev}(C^\infty(\Sigma_{g'})\rtimes G)$$ with $\C$ coefficients, via the
Chern character.

Thus orbifold line bundles on $\Sigma(g,\nu_1,\ldots,\nu_n)$ can be also
described as $G$-equivariant line bundles over the covering
smooth surface $\Sigma_{g'}$, and again as $\SO(2)$-equivariant line
bundles on $G\backslash \widehat P$.

\begin{rems} {\em
With the notation used in the previous section, let $G$ be a finite group
acting smoothly and properly on a smooth compact oriented manifold $X$.
There is a natural choice of a fundamental class $[X]_G \in H_0(X,G)$ in
the delocalized equivariant homology of $X$, given by the fundamental
classes of each compact oriented smooth submanifold $X^\gamma$, $[X]_G =
\oplus_{\gamma\in G} [X^\gamma]$. In the case of hyperbolic 2-orbifolds,
the equivariant fundamental class $[ \Sigma_{g'} ]_G$ is given by $$ [
\Sigma_{g'} ]_G=[\Sigma_{g'}]\oplus_{j} [p_j]^{\nu_j-1}\in
H_2(\Sigma_{g'},\C)\oplus_j (H_0(p_j,\C))^{\nu_j-1}.$$ The corresponding
equivariant Euler number $< ch_G({\cal L}), [\Sigma_{g'}]_G >$ is obtained
by evaluating $$ < ch_G({\cal L}), [\Sigma_{g'}]_G >= <c_1(\widetilde {\cal
L}), [\Sigma_{g'}]> + \sum_{j=1}^n \sum_{\gamma\in G_{p_j}\backslash \{1
\}} \lambda(\gamma). $$ }\end{rems}

\subsection{Classifying space of the orbifold fundamental group}

Here we find it convenient to follow
Baum, Connes and Higson \cite{BC},
\cite{BCH}. Let $M$ be a
good orbifold, that is its orbifold universal cover $\widetilde M$ is a
smooth manifold which has a
{\em proper} $\Gamma$-action, where $\Gamma$ denotes the orbifold
fundamental group of $M$. That is, the map $$
\begin{array}{rl}
\widetilde M \times \Gamma & \to \widetilde M \times \widetilde M\\ (x,
\gamma) & \to (x , \gamma x)
\end{array}
$$
is a proper map. The universal example for such a proper action is denoted
in \cite{BC}, \cite{BCH} by $\underline E\Gamma$. It is universal in the
sense that there is a continuous $\Gamma$-map $$ f: \widetilde M \to
\underline E\Gamma
$$ which is unique up to $\Gamma$-homotopy, and moreover $\underline
E\Gamma$ itself is unique up to $\Gamma$-homotopy. The quotient $\underline
B\Gamma = \Gamma\backslash \underline E\Gamma$ is an orbifold. Just as
$B\Gamma$ classifies isomorphism classes of $\Gamma$-covering spaces, it
can be shown that $ \underline B\Gamma$ classifies isomorphism classes of
orbifold $\Gamma$-covering spaces.

\begin{examples}{\em
It turns out that if $\Gamma$ is a discrete subgroup of a connected Lie
group $G$, then $\underline E\Gamma = G/K$, where $K$ is a maximal compact
subgroup.} \end{examples}

\begin{examples} {\em
The orbifold $\Sigma(g,\nu_1,\ldots,\nu_n)$, viewed as the quotient space
$$\Sigma(g,\nu_1,\ldots,\nu_n)=\Gamma\backslash {\bf H}$$ is an example
of the above construction.}
\end{examples}

This is the main class of examples that we are concerned with in this paper.

Let $S\Gamma$ denote the set of all elements of $\Gamma$ which are of
finite order. Then $S\Gamma$ is not empty, since $1\in S\Gamma$. $\Gamma$
acts on $S\Gamma$ by conjugation, and let $F\Gamma$ denote the associated
permutation module over $\C$, i.e. $$
F\Gamma = \left\{\sum_{\alpha \in S\Gamma} \lambda_\alpha [\alpha] \,
\Big|\, \lambda_\alpha \in \C \quad \hbox{and}\quad \lambda_\alpha = 0
\quad \hbox{except for
a finite number of}\quad\alpha\right\}
$$

\subsection{Twisting an elliptic operator} We will discuss elliptic
operators only on good orbifolds, and refer to \cite{Kaw} for the general
case.
Let $M$ be a good orbifold, that is the universal orbifold cover
$\widetilde M$ of $M$ is a smooth manifold. Let $\W\to \M$ be a
$\Gamma$-invariant Hermitian vector bundle over $\widetilde M$. Let $D$ be
a 1st order elliptic differential operator on $M$, $$ D: L^2 (M, {\cal E})
\to L^2(M, {\cal F}) $$ acting on $L^2$ orbifold sections of the orbifold
vector bundles ${\cal E}, {\cal F}$ over $M$. By fiat, $D$ is a
$\Gamma$-equivariant 1st order elliptic differential operator $\widetilde
D$ on the smooth manifold $\widetilde M$,
$$
\widetilde D :L^2 (\widetilde M, \E) \to L^2(\widetilde M, \F). $$ Given
any connection
$\nabla^{\W}$ on $\W$ which is compatible with the $\Gamma$ action and the
Hermitian metric, we wish to define an extension of the elliptic operator
$\widetilde D$, to act on sections of $\E\otimes \W$, $\F\otimes \W$. \[
\widetilde D\otimes\nabla^{\W}:\Gamma(M,\E\otimes \W)\to \Gamma(M,\F\otimes
\W) \]
and we want it to satisfy the following property:\ If $\W$ is a trivial
bundle, and $\nabla^0$ is the trivial connection on $\W$, then for $u \in
\Gamma(\M,\E),\
h \in \Gamma(\M, \W)$ such that $\nabla^0 h = 0$, \[ (\widetilde
D\otimes\nabla^0)(u\otimes h)= (\widetilde D u)\otimes h \]
To do this, define a morphism
$$
\begin{array}{rl}
&\quad S=S_D\,:\,\E\otimes T^* \M\to\F\\ &\quad S(u\otimes df) = \widetilde
D(f u)-f{\widetilde D}u \end{array} $$
for $f\in C^\infty(\M)$ and
$u\in\Gamma(M,\E)$. Then $S$ is a tensorial. Consider $S=S\otimes
1:\E\otimes T^* \M \otimes \W\to\F\otimes \W$ defined by \[ S(u\otimes
df\otimes e)=S(u\otimes df)\otimes e \] for $u, f$ as before and
$e\in\Gamma(M,\W)$.\\ Recall that a connection $\nabla^{\W}$ on $\W$ is a
derivation \[
\nabla^{\W}:\Gamma(\M,\W)\to\Gamma(\M,T^* \M\otimes \W) \] Define
$\widetilde D\otimes\nabla^{\W}$ as \[
(\widetilde D\otimes\nabla^{\W})(u\otimes e)=(D u)\otimes
e+S(u\otimes\nabla^{\W} e) \]
Then $\widetilde D\otimes\nabla^{\W}$ is a $1^{\hbox{ st }}$ order elliptic
operator.

\subsection{Twisted index theorem for orbifolds}

Let $M$ be a compact orbifold of dimension $n=4\ell$. Let
$\Gamma\to\widetilde{M}\stackrel{p}{\to}M$ be the universal orbifold cover
of $M$ and the orbifold fundamental group is $\Gamma$. Let $D$ be an
elliptic 1st order
operator on $M$, that is a $\widetilde D$ on $\widetilde M$, \[
\widetilde D\,:\,L^2(\widetilde M,
\E)\to
L^2(\widetilde M,\F), \]
such that $\widetilde D$ commutes with the $\Gamma$-action on $\widetilde{M}$.

Now let $\widetilde\omega$ be a $\Gamma$-invariant closed 2-form on
$\widetilde M$, $\widetilde\omega=d\eta$. Define $\nabla=d+\,i\eta$. Then
$\nabla$ is a Hermitian connection on the trivial line bundle over
$\widetilde{M}$, and the curvature of $\nabla,\ (\nabla)^2=i\,
\widetilde{\omega}$. (Here $s\in {\R}$.) Then $\nabla$ defines a projective
$(\Gamma,\sigma)$-action on $L^2$ spinors as in section 1.2.

Consider the {\em twisted elliptic operator} on $\widetilde{M}$, \[
\widetilde{D}\otimes\nabla\,:\,L^2(\widetilde{M}, \E)\to
L^2(\widetilde{M},\F) \]
Then $\widetilde{D}\otimes\nabla$ no longer commutes with $\Gamma$, but it
does commute with the projective $(\Gamma,\sigma)$ action. Let $P_+, P_-$
be the orthogonal projections onto the null space of
$\widetilde{D}\otimes\nabla$ and
$(\widetilde{D}\otimes\nabla)^*$ respectively since \[
(\widetilde{D}\otimes\nabla)\ P_+ = 0 \qquad {\hbox{and}} \qquad
(\widetilde{D}\otimes\nabla)^*\ P_- = 0\] By elliptic regularity, it
follows that the Schwartz (or integral) kernels of $P_\pm$ are smooth.
Since $\widetilde{D}\otimes\nabla$ and its adjoint commutes with the
$(\Gamma,\sigma)$ action, one has $$ e^{-i\phi_\gamma(x)} P_\pm(\gamma
x,\gamma y)\,e^{i\phi_\gamma(y)} = P_\pm(x,y)\quad\forall\gamma\in\Gamma. $$
In particular, $P_\pm(x,x)$ is smooth and $\Gamma$-invariant on
$\widetilde{M}$.
Therefore the corresponding von Neumann trace (cf. section 1.2) is {\em
finite}, \[
\tr\left(P_\pm\right) =
\int_M\,{\underline\tr} \left(P_\pm(x,x)\right)\,dx < \infty.\] The
$L^2$-index is by definition $$
\Index_{L^2} (\widetilde{D}\otimes\nabla) = \tr(P_+)-\tr(P_-). $$

To describe the next theorem, we will briefly review some material on
characteristic classes for orbifold vector bundles. Let $M$ be a good
orbifold, that is the universal orbifold cover $\Gamma \to \widetilde M \to
M$ of $M$
is a smooth manifold. Then the orbifold tangent bundle $TM$ of $M$, can be
viewed as the $\Gamma$-equivariant bundle $T\widetilde M$ on $\widetilde
M$. Similar comments apply to the orbifold cotangent bundle $T^*M$ and more
generally,
any orbifold vector bundle on $M$. It is then clear that choosing
$\Gamma$-invariant connections on the $\Gamma$-invariant vector bundles on
$\widetilde M$, one can define the Chern-Weil representatives of the
characteristic classes of the $\Gamma$-invariant vector bundles on
$\widetilde M$. These characteristic classes are $\Gamma$-invariant and so
define cohomology classes on $M$. For further details, see \cite{Kaw}.

\begin{thm}
Let $M$ be a compact, even dimensional, good orbifold, $\Gamma$ be its
orbifold fundamental group, $\widetilde D$ be a $\Gamma$-invariant twisted
Dirac operator on $\widetilde M$, where $\Gamma \to \widetilde M\to M$ is
the universal orbifold cover of $M$. Then one has
$$
index_{L^2}(\widetilde D\otimes \nabla) = \frac{q !} {(2\pi i)^q (2q!)}
\left<Td(M)\cup ch(symb(D))
\cup e^{\omega}, [T^*M] \right>
$$
where $Td(M)$ denotes the Todd characteristic class of the complexified
orbifold tangent bundle of $M$ which is pulled back to the orbifold
cotangent bundle $T^*M$, $ch(symb(D))$ is the Chern character of the symbol
of the operator $D$,
\label{thm1}
\end{thm}

\noindent{\bf Proof.}

The proof is similar to the case of
Atiyah's $L^2$ index theorem
for covering spaces. An important conceptual difference lies in the fact
that $\widetilde M$ is an orbifold cover, and not an actual cover of $M$.
We have
$\widetilde D = {\widetilde{\not\!\partial}^\pm \otimes\nabla^{\cal E}} =
\widetilde{\not\!\partial}^\pm_{\cal E} $. Let $k^\pm(t,x,y)$ denote the
heat kernel of the $\Gamma$-invariant Dirac operators
$(\widetilde{\not\!\partial}^\pm_{\cal E} \otimes \nabla)^2$ on the
universal orbifold cover of $M$, and $P^\pm(x,y)$ the smooth Schwartz
kernels o
f the orthogonal projections $P^\pm$ onto the null space of
$\widetilde{\not\!\partial}^\pm_{\cal E}\otimes \nabla^s$. By a general
result
of Cheeger-Gromov-Taylor \cite{CGT} (see also \cite{Roe}), the heat kernel
$k^\pm(t,x,y)$ converges uniformly over compact subsets of $\widetilde M
\times \widetilde M$ to $P^\pm(x,y)$, as $t \to \infty$. Therefore one has
\begin{equation}
\begin{array}{ll}
\lim_{t\to\infty} \tr(e^{-t(\widetilde{\not\partial}^\pm_{\cal E} \otimes
\nabla)^2})
& = \displaystyle\lim_{t\to\infty} \int_M \underline\tr(k^\pm(t,x,x))
dx\\[2mm]
& = \displaystyle\int_M \underline\tr(P^\pm(x,x)) dx\\[2mm]
& = \tr(P^\pm) \label{1}
\end{array}
\end{equation}

Next observe that
$$
\begin{array}{rl}
\frac{\partial}{\partial t} \tr_s( e^{-t (\widetilde{\not\partial}_{\cal E}
\otimes \nabla)^2})
& = - \tr_s( (\widetilde{\not\!\partial}_{\cal E} \otimes \nabla)^2
e^{-t(\widetilde{\not\partial}_{\cal E} \otimes \nabla)^2})\\ & = - \tr_s(
[\widetilde{\not\!\partial}_{\cal E} \otimes \nabla,
(\widetilde{\not\!\partial}_{\cal E} \otimes \nabla) e^{-t
(\widetilde{\not\partial}_{\cal E} \otimes \nabla)^2}])\\ & = 0
\end{array} $$
since $\widetilde{\not\!\partial}_{\cal E}\otimes \nabla$ is an odd
operator. Here $\tr_s$ denotes the graded trace, i.e. the composition of
the trace $\tr$ and the
grading operator. Therefore we deduce that \begin{equation}\begin{array}{ll}
\tr_s( e^{-t (\widetilde{\not\partial}_{\cal E} \otimes \nabla)^2}) & =
\lim_{t\to\infty} \tr_s( e^{-t (\widetilde{\not\partial}_{\cal E} \otimes
\nabla)^2})\\[2mm]
& = \tr_s (P) \\[2mm]
& = \Index_{L^2}(\widetilde{\not\!\partial}^+_{\cal E} \otimes \nabla).
\label{2}
\end{array}\end{equation}
By the local index theorem of Atiyah-Bott-Patodi \cite{ABP}, Getzler
\cite{Get}, one has
\begin{equation}
\lim_{t\to 0}\left( {\underline\tr}(k^+(t,x,x)) -
{\underline\tr}(k^-(t,x,x)) \right) =
[\widehat{A}(\Omega)\,\tr(e^{R^{\E}}) e^{\omega}]_n \label{(3)}
\end{equation}
where $[\;\;]_n$ denotes the component of degree $n=$ dim $M$, $\Omega$ is
the curvature of the metric on $\widetilde M$, $R^{\E}$ is the curvature of
the connection on $\widetilde {\cal E}$. Combining equations $(1.3), (1.4)$
and $(1.5)$, one has
$$
\Index_{L^2}(\widetilde{\not\!\partial}^+_{\cal E}\otimes \nabla) =
\displaystyle\int_M \widehat{A}(\Omega)\,\tr(e^{R^{\E}}) e^{\omega}. $$

$\diamond$

\begin{rems}{\em
A particular case of Theorem \ref{thm1} highlights a key new phenomenon in
the case of orbifolds, viz. in the special case when the multiplier
$\sigma=1$ is trivial, then
the $\Index_{L^2}(\widetilde D)$ formally coincides the $L^2$ index of
$\widetilde D$ as defined Atiyah \cite{At}. By comparing with the
cohomological formula due to Kawasaki \cite{Kaw} for the Fredholm index of
the operator $ D$ on the orbifold $M$, we see that in general these are
{\em not} equal, and the error term is a rational number which can be
expressed explicitly as a cohomological formula on the lower dimensional
strata of the orbifold $M$. Thus, we see that for general orbifolds the
$L^2$ index of $\widetilde D$, is only a {\em rational} number. This was
also observed by \cite{Far}. This is in contrast to the situation
when the orbifold is smooth, where  Atiyah's $L^2$ index theorem
establishes the integrality of the $L^2$ index in this case.}
\end{rems}

\section{Range of the trace and the Kadison constant}

In this section, we will first calculate the range of the canonical trace
map on $K_0$ of the twisted group $C^*$-algebras for Fuchsian groups
$\Gamma$ of signature $(g,\nu_1,\ldots ,\nu_n)$. We use in an essential way
some of the results of the previous section such as the twisted version of
the $L^2$-index theorem of Atiyah \cite{At}, which is due to Gromov
\cite{Gr2}, and which is proved in Theorem \ref{thm1}. This enables us to
deduce information about projections in the twisted group $C^*$-algebras.
In the case of no twisting, this follows because the Baum-Connes conjecture
is known to be true while these results are also well known for the case of
the irrational rotation algebras, and for the twisted groups $C^*$ algebras
of the fundamental groups of closed Riemann surfaces of positive genus
\cite{CHMM}. Our theorem generalizes most of these results. Moreover, we
prove analogous results in the case of compact, real and complex hyperbolic
four-manifolds.
We will apply the results of this section in the next section to study some
quantitative aspects of the
spectrum of projectively periodic elliptic operators, mainly on orbifold
covering spaces
of hyperbolic orbifolds.

\subsection{The isomorphism classes of algebras ${C}^*(\Gamma, \sigma)$}
Let $\sigma \in Z^2(\Gamma, {\mathbf U}(1))$ be a multiplier on $\Gamma$,
where $\Gamma$ is a Fuchsian group of signature $(g, \nu_1, \ldots,
\nu_n)$. If $\sigma' \in Z^2(\Gamma, {\mathbf U}(1))$ is another multiplier
on $\Gamma$ such that $[\sigma] = [\sigma'] \in H^2(\Gamma, {\mathbf
U}(1))$, then it can be easily shown that ${C}^*(\Gamma, \sigma) \cong
{C}^*(\Gamma, \sigma')$. That is, the isomorphism classes of the
$C^*$-algebras ${C}^*(\Gamma_{g'}, \sigma)$ are naturally parameterized by
$H^2(\Gamma, {\mathbf U}(1))$. In particular, if we consider only
multipliers $\sigma$ such that $\delta(\sigma) = 0$, we see that these are
parameterized by $\hbox{ker} (\delta) \subset H^2(\Gamma, {\mathbf U}(1))$.
It follows from the discussion at the beginning of the next subsection that
$\hbox{ker} (\delta) \cong {\mathbf U}(1)$. We summarize this below.

\begin{lemma} Let $\Gamma$ be a Fuchsian group of signature $(g, \nu_1,
\ldots, \nu_n)$.
Then the isomorphism classes of twisted group $C^*$-algebras ${C}^*(\Gamma,
\sigma)$ such that $\delta(\sigma) = 0$ are naturally parameterized by
${\mathbf U}(1)$.
\end{lemma}

\subsection{K-theory of twisted group $C^*$ algebras} We begin by computing
the $K$-theory of twisted group $C^*$-algebras for Fuchsian groups $\Gamma$
of signature $(g,\nu_1,\ldots ,\nu_n)$. Let $\sigma$ be a multiplier on
$\Gamma$. It defines a cohomology class $[\sigma] \in H^2(\Gamma, {\mathbf
U}(1))$. Consider now the short
exact sequence of coefficient groups
\[
1\to {\Z} \stackrel{i}{\to} {\R}
\stackrel{e^{2\pi\sqrt{-1}}}{\longrightarrow} {\mathbf U}(1) \to 1, \]
which gives rise to a long exact sequence of cohomology groups (the change
of coefficient groups sequence)

\begin{equation}
\cdots \to H^2(\Gamma,{\Z}) \stackrel{i_*}{\to} H^2(\Gamma,{\R})
\stackrel{{e^{2\pi\sqrt{-1}}}_*}{\longrightarrow} H^2(\Gamma, {\mathbf
U}(1)) \stackrel{\delta}{\to} H^3(\Gamma,{\Z})\stackrel{i_*}{\to}
H^3(\Gamma,{\R}).
\label{(3.1)}
\end{equation}

We first show that the map
$$ H^2(\Gamma, {\mathbf U}(1))\stackrel{\delta}{\to} H^3(\Gamma,{\Z})$$ is
a a surjection.

In fact, it is enough to show that $H^3(\Gamma,{\R})=\{0\}$. In order to
see this it is enough to notice that we have a $G$ action on $B\Gamma_{g'}$
with quotient $B\Gamma$, \begin{equation} G\to B\Gamma_{g'}
\stackrel{\lambda}{\to} B\Gamma \label{(3.2)} \end{equation}
and therefore, in the Leray-Serre spectral sequence, we have $$E^2=Tor^{H_*(G,
\R)}(\R,H_*(B\Gamma_{g'},\R)) $$ that converges to $H_*(B\Gamma,\R)$.
Moreover, we have $$ E^2=Tor^{H_*(G, \R)}(\R,\R) $$ converging to
$H_*(BG,\R)$, see 7.16 of \cite{McCl}.

Notice also that, with $\R$ coefficients,
we have $ H_q(BG,\R)=\{0\}$ for $q>0$. Thus we obtain that, with $\R$
coefficients, $H_q(B\Gamma,\R)\cong H^q(B\Gamma, \R)$
is $\R$ in degrees $q=0$ and $q=2$, $\R^{2g}$ in degree $q=1$, and trivial
in degrees $q>2$. In particular, (\ref{(3.1)}) now becomes \begin{equation}
\cdots \to H^2(\Gamma,{\Z}) \stackrel{i_*}{\to} H^2(\Gamma,{\R})
\stackrel{{e^{2\pi\sqrt{-1}}}_*}{\longrightarrow} H^2(\Gamma, {\mathbf
U}(1)) \stackrel{\delta}{\to} H^3(\Gamma,{\Z})\stackrel{i_*}{\to}
0.\label{(3.3)}
\end{equation}

In the following $[\Gamma]$ will denote a choice of a generator in
$$H_2(B\Gamma,\R)\cong\R\cong H^2( B\Gamma, \R).$$  Using equation
(\ref{(3.2)}) and the previous argument, we see that
$\lambda_*[\Sigma_{g'}] = \#(G) [\Gamma]$, since $B\Gamma_{g'}$ and
$\Sigma_{g'}$ are homotopy equivalent,
and where $ \#(G)$ denotes the order of the finite group $G$.

In particular, for any multiplier $\sigma$ of $\Gamma$ with $ [\sigma] \in
H^2(\Gamma_{g'}, {\mathbf U}(1))$ and with $\delta(\sigma) =0$, there is a
$\R$-valued 2-cocycle $\zeta$ on $\Gamma$ with $[\zeta] \in H^2(\Gamma,
\R)$ such that
$[e^{2\pi\sqrt{-1}\zeta}] = [\sigma]$. Define a homotopy $[\sigma_t] =
[e^{2\pi\sqrt{-1}t\zeta}] \quad \forall t \in [0,1]$ which is a homotopy of
multipliers $\sigma_t$ that connects the multiplier $\sigma$ and the
trivial multiplier. Note also that this homotopy is {\em canonical} and not
dependent on the particular choice of $\zeta$. Therefore one obtains a
homotopy of isomorphism classes of twisted group $C^*$-algebras
$C^*(\Gamma, \sigma_t)$ connecting $C^*(\Gamma, \sigma)$ and $C^*(\Gamma)$.
It is this homotopy which will essentially be used to show that
$C^*(\Gamma, \sigma)$ and $C^*(\Gamma)$ have the same $K$-theory.

Let $\Gamma\subset G$ be a discrete
cocompact subgroup of $G$ and $A$ be an algebra admitting an action of
$\Gamma$ by automorphisms. Then the cross product algebra $[A\otimes
C_0(G)]\rtimes \Gamma$, is Morita equivalent to the algebra of continuous
sections vanishing at infinity
$C_0(\Gamma\backslash G, {\cal E})$, where ${\cal E}\to \Gamma\backslash G$
is the flat $A$-bundle defined as the quotient \begin{equation}
{\cal E} = (A\times G)/\Gamma \to \Gamma\backslash G. \end{equation} Here
we consider the diagonal action of $\Gamma$ on $A\times G$. We refer the
reader to \cite{Kas1} for the technical definition of a $K$-amenable group.
However we mention that any solvable Lie group, and in fact any amenable
Lie group is $K$-amenable, and in fact it is shown in \cite{Kas1},
\cite{JuKas} that the
non-amenable groups ${\mathbf{\SO}}_0(n,1)$, ${\mathbf{SU}}(n,1)$ are
$K$-amenable Lie groups. Also,
Cuntz \cite{Cu} has shown that the class of $K$-amenable groups is closed
under the operations of taking subgroups, under free products and under
direct products.

\begin{thm}[\cite{Kas1},\cite{Kas2}]
If $G$ is $K$-amenable, then
$(A\rtimes\Gamma)\otimes C_0(G)$ and $[A\otimes C_0(G)]\rtimes \Gamma$ have
the same $K$-equivariant $K$-theory, where $K$ acts in the standard way on
$G$ and trivially on the other factors. \label{thm2} \end{thm}

Combining Theorem \ref{thm2} with the remarks above, one gets the following
important corollary.

\begin{cor}
If $G$ is $K$-amenable, then $(A\rtimes\Gamma)\otimes C_0(G)$ and
$C_0(\Gamma\backslash G, {\cal E})$ have the same $K$-equivariant
$K$-theory.
Equivalently, one has for $j=0,1$,
$$
{K_K}_j(C_0(\Gamma\backslash G, {\cal E})) \cong {K_K}_{j+ \dim(G/K)}
(A\rtimes\Gamma).
$$
\label{cor3}
\end{cor}

We now come to the main theorem of this section, which generalizes theorems
of \cite{CHMM}, \cite{PR}, \cite{PR2}.

\begin{thm}
\label{thm3}
Suppose given $\Gamma$  a discrete cocompact subgroup in a $K$-amenable
Lie group $G$ and suppose that $K$ is a maximal compact subgroup of
$G$. Then $$ {K}_\bullet(C^*(\Gamma,\sigma)) \cong
{K_K}^{\bullet + \dim(G/K)}(\Gamma\backslash G, \delta(B_\sigma)), $$ where
$\sigma \in H^2(\Gamma, {\mathbf U}(1))$ is any multiplier on
$\Gamma$. Here
${K_K}^{\bullet}(\Gamma\backslash G,
\delta(B_\sigma))$ is the twisted $K$-equivariant $K$-theory of a
continuous trace
$C^*$-algebra
$B_\sigma$ with spectrum $\Gamma\backslash G$, and $\delta(B_\sigma)$
denotes the
Dixmier-Douady invariant of $B_\sigma$.\end{thm}

\noindent{\bf Proof.}
Let $\sigma\in H^2(\Gamma, {\mathbf U}(1))$, then the {\em twisted} cross
product algebra $A\rtimes_\sigma \Gamma$ is stably equivalent to the cross
product $(A\otimes
{\cal K})\rtimes\Gamma$ where ${\cal K}$ denotes compact operators. This is
the Packer-Raeburn stabilization trick \cite{PR}, which we now describe in
more detail.
Let $V:
\Gamma \to U(\ell^2(\Gamma))$ denote the left regular $(\Gamma,
\bar\sigma)$ representation on $\ell^2(\Gamma)$, i.e. for $\gamma, \gamma_1
\in \Gamma$ and $f\in \ell^2(\Gamma)$
$$
(V(\gamma_1)f)(\gamma) = \bar\sigma(\gamma_1, \gamma_1^{-1}\gamma)
f(\gamma_1^{-1}\gamma).
$$
Then for $\gamma_1, \gamma_2 \in \ell^2(\Gamma)$, $V$ satisfies
$V(\gamma_1)V(\gamma_2) = \bar\sigma(\gamma_1,
\gamma_2)V(\gamma_1\gamma_2)$. That is, $V$ is a projective representation
of $\Gamma$. Since $Ad$ is trivial on ${\mathbf U}(1)$, it follows that
$\alpha(\gamma) = Ad(V(\gamma))$ is a representation of $\Gamma$ on ${\cal
K}$. This is easily generalized to the case when ${\C}$ is replaced by the
$*$ algebra $A$.

Using Corollary \ref{cor3},
one sees
that $A\rtimes_\sigma\Gamma \otimes C_0(G)$ and $C_0(\Gamma\backslash G,
{\cal E}_\sigma)$
have the same $K$-equivariant $K$-theory, whenever $G$ is $K$-amenable,
where $$
{\cal E}_\sigma = (A\otimes{\cal K} \times G)/\Gamma \to \Gamma\backslash G
$$
is a flat $A\otimes{\cal K}$-bundle over $\Gamma\backslash G$ and $K$ is a
maximal compact subgroup of $G$. In the particular case when $A= {\C}$, one
sees that $C^*_r(\Gamma, \sigma)\otimes C_0(G)$ and $C_0(\Gamma\backslash
G, {\cal E}_\sigma)$ have the same $K$-equivariant $K$-theory whenever $G$
is $K$-amenable, where \[ {\cal E}_\sigma = ({\cal K} \times
G)/\Gamma\to\Gamma\backslash G. \] But the twisted $K$-equivariant
$K$-theory ${K_K}^{*}(\Gamma\backslash G, \delta(B_\sigma))$ is by
definition the same as the $K$-equivariant $K$-theory of the continuous trace
$C^*$-algebra $B_\sigma = C_0(\Gamma\backslash G, {\cal E}_\sigma)$ with
spectrum $\Gamma\backslash G$.
Therefore
$$
K_\bullet(C^*(\Gamma,\sigma)) \cong {K_K}^{\bullet + \dim(G/K)}
(\Gamma\backslash G, \delta(B_\sigma)).
$$

$\diamond$

Our next main result says that for discrete cocompact subgroups in
$K$-amenable Lie groups,
the reduced and unreduced twisted group $C^*$-algebras have canonically
isomorphic $K$-theories. Therefore all the results that we prove regarding
the $K$-theory of these reduced twisted group $C^*$-algebras are also valid
for the unreduced twisted group $C^*$-algebras.

\begin{thm}
Let $\sigma\in H^2(\Gamma, {\mathbf U}(1))$ be a multiplier on $\Gamma$ and
$\Gamma$ be a discrete cocompact
subgroup in a $K$-amenable Lie group. Then the canonical morphism
$C^*(\Gamma, \sigma) \rightarrow C^*_r(\Gamma, \sigma)$ induces an
isomorphism
$$
K_*(C^*(\Gamma, \sigma)) \cong K_*(C^*_r(\Gamma, \sigma)). $$ \label{thm4}
\end{thm}

\noindent{\bf Proof.}
We note that by the Packer-Raeburn trick, one has $$ C^*(\Gamma, \sigma)
\otimes {\cal K} \cong {\cal K} \rtimes \Gamma $$ and
$$
C^*_r(\Gamma, \sigma) \otimes {\cal K} \cong {\cal K} \rtimes_r \Gamma, $$
where $\rtimes_r$ denotes the reduced crossed product. Since $\Gamma$ is a
lattice in a $K$-amenable Lie group, the canonical morphism ${\cal K}
\rtimes \Gamma \rightarrow {\cal K} \rtimes_r \Gamma$ induces an
isomorphism (cf. \cite{Cu})
$$
K_*({\cal K} \rtimes \Gamma) \cong K_*({\cal K} \rtimes_r \Gamma), $$ which
proves the result.
$\diamond$

We now specialize to the case where we have  $G= \so_0(2,1)$,
$K=\so(2)$ and $\Gamma
= \Gamma(g, \nu_1, \ldots, \nu_n)$ is a Fuchsian group, i.e. the orbifold
fundamental group of a hyperbolic orbifold of signature $(g, \nu_1, \ldots,
\nu_n)$, \ $\Sigma(g, \nu_1, \ldots, \nu_n)$, where $\Gamma \subset G$
(note that $G$ is $K$-amenable), or when $G= {\R}^2$, $K= \{e\}$ and $g=1$,
with $\Gamma$ being a cocompact crystallographic group.

\begin{prop}
Let $\sigma$ be a multiplier on the Fuchsian group $\Gamma$ of signature
$(g,\nu_1,\ldots ,\nu_n)$ such that $\delta(\sigma) = 0$. Then one has
\begin{enumerate}
\item $K_0(C^*(\Gamma, \sigma)) \cong K_0(C^*(\Gamma)) \cong
K^0_{orb}(\Sigma(g, \nu_1,\ldots,\nu_n)) \cong \Z^{2-n +\sum_{j=1}^n
\nu_j}$ \item $K_1(C^*(\Gamma, \sigma)) \cong K_1 (C^*(\Gamma)) \cong
K^1_{orb}(\Sigma(g, \nu_1,\ldots,\nu_n)) \cong {\Z}^{2g}$. \end{enumerate}
\label{prop4}
\end{prop}

\noindent{\bf Proof.}
Now by a result due to Kasparov \cite{Kas1}, which he proves by connecting
the regular representation to the trivial one via the complementary series,
one has
$$
K_\bullet (C^*(\Gamma)) \cong K^\bullet_{\SO(2)}(P(g, \nu_1,\ldots,\nu_n))
= K^\bullet_{orb}(\Sigma(g,\nu_1,\ldots ,\nu_n)). $$ We recall next the
calculation of Farsi \cite{Far} for the orbifold $K$-theory of the
hyperbolic orbifold $\Sigma(g,\nu_1,\ldots ,\nu_n)$ $$ \begin{array}{ll}
K^0_{orb}(\Sigma(g, \nu_1,\ldots,\nu_n)) & \equiv K_0(C^*(\Sigma(g,
\nu_1,\ldots,\nu_n)) ) \\
& =
K^0_{\SO(2)}(P(g, \nu_1,\ldots,\nu_n))
\cong \Z^{2-n +\sum_{j=1}^n \nu_j} \end{array} $$
and
$$ \begin{array}{ll}
K^1_{orb}(\Sigma(g, \nu_1,\ldots,\nu_n)) & \equiv K_1(C^*(\Sigma(g,
\nu_1,\ldots,\nu_n)) ) \\ & =
K^1_{\SO(2)}(P(g, \nu_1,\ldots,\nu_n))
\cong \Z^{2g}\end{array}
$$

By Theorem \ref{thm3} we have
\[
K_j(C^*(\Gamma)) \cong K^j_{\SO(2)}(P(g, \nu_1,\ldots,\nu_n))
\quad\hbox{for }j=0,1,
\]
and more generally
\[
K_j(C^*(\Gamma,\sigma)) \cong K^j_{\SO(2)}(P(g, \nu_1,\ldots,\nu_n),
\delta(B_\sigma)),\quad
j = 0,1,
\]
where $B_\sigma = C(P(g, \nu_1,\ldots,\nu_n), {\cal E}_\sigma)$. Finally,
because ${\cal E}_\sigma$ is a locally trivial bundle of $C^*$-algebras
over $P(g, \nu_1,\ldots,\nu_n)$, with fiber ${\cal K}$ ($=$ compact
operators), it has a Dixmier-Douady invariant $\delta(B_\sigma)$ which can
be viewed as the obstruction to $B_\sigma$ being Morita equivalent to
$C(\Sigma_g)$. But by assumption $\delta(B_\sigma) = \delta(\sigma) = 0$.
Therefore $B_\sigma$ is Morita equivalent to $C(P(g, \nu_1,\ldots,\nu_n))$
and we conclude that \[
K_j(C^*(\Gamma, \sigma)) \cong K^j_{\SO(2)}(P(g, \nu_1,\ldots,\nu_n)) \cong
K^j_{orb}(\Sigma(g, \nu_1,\ldots,\nu_n))\quad j = 0,1. \] $\diamond$

\subsection{Twisted Kasparov map}
Let $\Gamma$ be as before,
that is, $\Gamma$ is the orbifold fundamental group of the hyperbolic
orbifold $\Sigma(g, \nu_1,\ldots,\nu_n)$. Then for any multiplier $\sigma$
on $\Gamma$, the \emph{twisted Kasparov isomorphism}, \begin{equation}
\mu_\sigma : K^\bullet_{orb} (\Sigma(g, \nu_1,\ldots,\nu_n)) \to K_\bullet
(C^*_r(\Gamma,\sigma))\label{(3.4)} \end{equation} is defined as follows.
Let
$${\cal E}\to\Sigma(g, \nu_1,\ldots,\nu_n)$$ be an orbifold vector bundle
over $\Sigma(g, \nu_1,\ldots,\nu_n)$ defining an element $$[{\cal E}] \in
K^0(\Sigma(g, \nu_1,\ldots,\nu_n)).$$ As in \cite{Kaw}, one can form the
twisted Dirac
operator
$$\npartial^+_{{\cal E}} : L^2(\Sigma(g, \nu_1,\ldots,\nu_n), {\cal
S}^+\otimes{\cal E})
\to L^2(\Sigma(g, \nu_1,\ldots,\nu_n), {\cal S}^-\otimes{\cal E})$$ where
${\cal S}^{\pm}$ denote the $\frac12$ spinor bundles over $\Sigma(g,
\nu_1,\ldots,\nu_n)$. By Proposition \ref{prop4} of the previous
subsection, there is a canonical isomorphism $$
K_\bullet (C^*_r(\Gamma, \sigma)) \cong K^\bullet_{orb} (\Sigma(g,
\nu_1,\ldots,\nu_n)).
$$
Both of these maps are assembled to yield the twisted Kasparov map as in
(\ref{(3.4)}).
Observe that $\Sigma(g, \nu_1,\ldots,\nu_n) = \underline B\Gamma$, and that
the twisted Kasparov map has a natural generalization, which will be
studied elsewhere.

We next describe this map more explicitly. One can lift the twisted Dirac
operator $\npartial^+_{{\cal E}}$ as above, to a $\Gamma$-invariant
operator $\widetilde{\npartial_{{\cal E}}^+}$ on ${\bf H} =
\widetilde\Sigma(g, \nu_1,\ldots,\nu_n)$, which is the universal orbifold
cover of $\Sigma(g, \nu_1,\ldots,\nu_n)$, \[ \widetilde{\npartial_{{\cal
E}}^+} : L^2({\bf H}, \widetilde{{\cal S}^+\otimes{\cal E}}) \to L^2({\bf
H}, \widetilde{{\cal S}^-\otimes{\cal E}}) \]
Therefore as before in (\ref{(3.3)}),
for any multiplier $\sigma$ of $\Gamma$ with $\delta([\sigma]) =1$, there
is a $\R$-valued 2-cocycle $\zeta$ on $\Gamma$ with $[\zeta] \in
H^2(\Gamma, \R)$ such that
$[e^{2\pi\sqrt{-1}\zeta}] = [\sigma]$. By the earlier argument using
spectral sequences and the fibration as in equation (\ref{(3.2)}), we see
that the map $ \lambda$ induces an isomorphism $H^2(\Gamma, \R) \cong
H^2(\Gamma_{g'}, \R)$, and therefore there is a 2-form $\omega$ on
$\Sigma_{g'}$ such that $[e^{2\pi\sqrt{-1}\omega}] = [\sigma]$.
Of course, the choice of $\omega$ is not unique, but this will not affect
the results that we are concerned with. Let $\widetilde \omega$ denote the
lift of $\omega$ to the universal cover ${\bf H}$. Since the hyperbolic
plane ${\bf H}$ is contractible, it follows that $\widetilde \omega =
d\eta$ where $\eta$ is a 1-form on ${\bf H}$ which is not in general
$\Gamma$ invariant. Now $\nabla = d - i\eta$ is a Hermitian connection on
the trivial complex line bundle on ${\bf H}$. Note that the curvature of
$\nabla$ is $\nabla^2 = i\widetilde\omega$. Consider now the twisted Dirac
operator $\widetilde \npartial^+_{{\cal E}}$ which is twisted again by the
connection $\nabla$, \[
\widetilde{\npartial_{{\cal E}}^+}\otimes\nabla : L^2({\bf H},
\widetilde{{\cal S}^+\otimes{\cal E}}) \to L^2({\bf H}, \widetilde{{\cal
S}^-\otimes{\cal E}}). \]
It does not commute with the $\Gamma$ action, but it does commute with the
projective $(\Gamma, \sigma)$-action which is defined by the multiplier
$\sigma$, and it has an
$(\Gamma, \sigma)$-$L^2$-index
\[
{\rm ind}_{(\Gamma, \sigma)} (\widetilde{\npartial^+_{{\cal E}}}\otimes
\nabla) \in
K_0(C^*_r(\Gamma, \sigma)).
\]
Formally, ${\rm ind}_{(\Gamma, \sigma)}
(\widetilde{\npartial^+_{{\cal E}}}\otimes \nabla) = [P^+] - [P^-]$, where
$P^\pm$ denotes the projection to the $L^2$ kernel of
$\widetilde{\npartial^\pm_{{\cal E}}}\otimes \nabla$. The problem is that
in general,
$P^\pm$ only lies in the twisted von Neumann algebra, and therefore one has
to add
a compact perturbation in $C^*_r(\Gamma, \sigma)$ to the operator, in order
to properly define the index. This is essentially the $C^*$ index of
Mishchenko-Fomenko
\cite{MiFo}, see also \cite{Ros}.
Then observe that the {\em twisted Kasparov map} is $$ \mu_\sigma([{\cal
E}]) ={\rm ind}_{(\Gamma, \sigma)} (\widetilde{\npartial^+_{{\cal
E}}}\otimes \nabla) ) \in K_0({C}^*(\Gamma, \sigma)).
$$

The canonical trace on ${C}^*_r(\Gamma, \sigma))$ induces a linear map $$
[\tr] : K_0 ({C}^*_r(\Gamma, \sigma)) \to {\R} $$ which is called the {\em
trace map} in $K$-theory. Explicitly, we first extend $\tr$ to matrices
with entries in ${C}^*(\Gamma, \sigma)$ as (with Trace denoting matrix
trace): $\tr(f\otimes r) = {\mbox{Trace}}(r) \tr(f). $
Then the extension of $\tr$ to $K_0$ is given by $[\tr]([e]-[f]) = \tr(e) -
\tr(f)$, where $e,f$ are idempotent matrices with entries in ${C}^*(\Gamma,
\sigma))$.

Clearly one has
$$
\Index_{L^2} (\widetilde{\npartial^+_{{\cal E}}}\otimes \nabla) ) =
[\tr]\left({\rm ind}_{(\Gamma, \sigma)} (\widetilde{\npartial^+_{{\cal
E}}}\otimes \nabla) )
\right)
$$
\subsection{Range of the trace map on $K_0$: the case of two dimensional
orbifolds}
We can now state the first major theorem of this section.

\begin{thm} Let $\Gamma$ be a Fuchsian group of signature $(g, \nu_1,
\ldots , \nu_n)$, and $\sigma$ be a multiplier of $\Gamma$ such that
$\delta(\sigma) = 0$.
Then the range of the trace map is
$$
[\tr] (K_0 ({C}^*_r(\Gamma, \sigma)) ) = \Z\theta + \Z + \sum_{i=1}^n \Z
(1/\nu_i) ,
$$
where $2\pi\theta = \langle[\sigma], [\Gamma]\rangle\ \in (0,1]$ is the
result of
pairing the
multiplier $\sigma$ with the fundamental class of $\Gamma$ (cf. subsection
2.1).
\label{thm5}
\end{thm}

\noindent{\bf Proof.}

We first observe that by the results of the previous subsection the twisted
Kasparov map is an isomorphism. Therefore to compute the range of the trace
map on $K_0$, it suffices to compute the range of the trace map on elements
of the form
$$\mu_\sigma([{\cal E}^0 ] -
[{\cal E}^1])$$
for any element
$$[{\cal E}^0 ] -
[{\cal E}^1] \in K^0_{orb}(\Sigma(g, \nu_1,\ldots,\nu_n)).$$ where ${\cal
E}^0, {\cal E}^1$ are orbifold vector bundles over the orbifold $\Sigma(g,
\nu_1,\ldots,\nu_n)$, which as in section 1.4, can be viewed
as $G$-equivariant vector bundles over the Riemann surface $\Sigma_{g'}$
which is
an orbifold $G$ covering of the orbifold $\Sigma(g, \nu_1,\ldots,\nu_n)$.

By the twisted $L^2$ index theorem for orbifolds, Theorem 1.7, one has
\begin{equation}
[\tr]({\rm ind}_{(\Gamma, \sigma)} (\widetilde{\npartial_{{\cal
E}}^+}\otimes \nabla)) =
\frac{1}{2\pi}\int_{\Sigma(g, \nu_1,\ldots,\nu_n)} \widehat{A}(\Omega)
\tr(e^{R^{{\cal E}}}) e^{\omega}.
\label{(3.5)}
\end{equation}
We next simplify the right hand side of equation (\ref{(3.5)}) using $$
\begin{array}{rl}
\widehat{A}(\Omega) &= 1 \\
\tr(e^{R^{{\cal E}}}) &= \rank {\cal E} + \tr(R^{{\cal E}}) \\ e^{\omega}
&= 1 + {\omega} .
\end{array} $$
Therefore one has
\[
[\tr]({\rm ind}_{(\Gamma, \sigma)} (\widetilde{\npartial_{{\cal
E}}^+}\otimes \nabla)) =
\frac{\rank {\cal E}}{2\pi}\;{\displaystyle \int_{\Sigma(g,
\nu_1,\ldots,\nu_n)} \omega} + \frac{1}{2\pi} { \displaystyle
\int_{\Sigma(g, \nu_1,\ldots,\nu_n)}\tr(R^{{\cal E}})}, \] Now by the index
theorem for orbifolds, due to Kawasaki \cite{Kaw}, we see that $$
\frac{1}{2\pi} { \displaystyle \int_{\Sigma(g,
\nu_1,\ldots,\nu_n)}\tr(R^{{\cal E}})} + \frac{1}{2\pi} {\sum_{i=1}^n
\beta_i/\nu_i } = \hbox{index} (\npartial_{{\cal E}}^+)
\in \Z,
$$
Therefore we see that
$$
\frac{1}{2\pi}{\displaystyle \int_{\Sigma(g,
\nu_1,\ldots,\nu_n)}\tr(R^{{\cal E}})} \in \Z + \sum_{i=1}^n \Z (1/\nu_i)
$$
Observe that
$$
\int_{\Sigma(g, \nu_1,\ldots,\nu_n)} \omega = \frac{1}{\#(G)
}\int_{\Sigma_{g'}} \omega = \langle [\omega], [\Gamma] \rangle $$ since
$\Sigma_{g'}$ is
an orbifold $G$ covering of the orbifold $\Sigma(g, \nu_1,\ldots,\nu_n)$
and $[\Gamma]$ is equal to $\frac{[\Sigma_{g'}]}{\#(G)}$, cf. section 2.1
and that by assumption,
$$\frac{\langle [\omega], [\Gamma] \rangle}{2\pi} - \theta\in {\Z}.$$ It
follows that the range of the trace map on $K_0$ is $${\Z} \frac{\langle
[\omega],
[\Gamma] \rangle}{2\pi} + \Z + \sum_{i=1}^n \Z (1/\nu_i)= {\Z}\theta + \Z +
\sum_{i=1}^n \Z (1/\nu_i).$$ $\diamond$

We will now discuss one application of this result, leaving further
applications to the
next section.The application studies the number of projections in the
twisted group $C^*$-algebra, which is a problem of independent interest.

We first recall the definition of the Kadison constant of a twisted group
$C^*$-algebra.
The {\em Kadison constant} of ${C}^*_r(\Gamma, \sigma)$ is defined by: $$
C_\sigma(\Gamma) = \inf\{ \tr(P) : P \ \ {\mbox{is a non-zero projection
in}} \ \
{C}^*_r(\Gamma, \sigma) \otimes {\cal K}\}. $$

\begin{prop} Let $\Gamma$ be as in Theorem 2.15. Let $\sigma$ be a
multiplier on
$\Gamma$ such that $\delta(\sigma) = 0$, and $2\pi\theta= \langle[\sigma],
[\Gamma]\rangle \in (0,1]$ be the result of pairing the cohomology class of
$\sigma$ with the fundamental class of $\Gamma$. If $\theta$ is rational,
then
there are at most a finite number of unitary equivalence classes of
projections, other than $0$ and $1$,
in the reduced twisted group $C^*$-algebra ${C}^*_r(\Gamma, \sigma)$.
\label{prop5}
\end{prop}

\noindent{\bf Proof.} By assumption, $\theta = p/q$. Let $P$ be a
projection in ${C}^*_r(\Gamma, \sigma)$. Then $1-P$ is also a projection in
${C}^*_r(\Gamma, \sigma)$ and one has $$ 1= \tr(1) = \tr(P) + \tr(1-P).
$$
Each term in the above equation is non-negative. Since $\sigma$ is rational
and by
Theorem \ref{thm5}, it follows that the Kadison constant
$C_\sigma(\Gamma)>0$ and $\tr(P) \in \{0, C_\sigma(\Gamma),
2C_\sigma(\Gamma), \ldots 1\}$. By faithfulness and normality of the
trace $\tr$, it follows that there are at most a finite number of unitary
equivalence classes of projections, other than those of $0$ and $1$ in
${C}^*_r(\Gamma, \sigma)$.
$\diamond$

\subsection{Range of the trace map on $K_0$: the case of 4 dimensional real
and complex hyperbolic manifolds}

We prove here the analogue of the results of the previous subsection, for
the case of compact, real and complex hyperbolic four dimensional manifolds.

We now set some notation for the theorem below. Let $\Gamma$ be a discrete,
torsion-free
cocompact subgroup of $G= \so_0(1, 4)$ or $\SU(1,2)$. We will assume that
$\delta(\sigma) = 0$, therefore there is a closed
$\Gamma$-invariant two form $\omega$ on $G/K$, where $K$ is a maximal
compact subgroup of $G$, such that $[e^{2i\pi \omega}] = [\sigma]$. Let
$Q(a,b) = \langle a\cup b, [\Gamma]\rangle\quad a,b \in H^2(\Gamma, \R)$ be
the intersection form on $\Gamma\backslash G/K$. Define the linear
functional $T_\omega : H^2(\Gamma, \Z) \to \R$ as $T_\omega (a) = Q(\omega,
a)$. Then we have:

\begin{thm} Let $\Gamma$ be a discrete, torsion-free cocompact subgroup of
$\so_0(1, 4)$ or of $\SU(1,2)$, and $\sigma$ be a multiplier of $\Gamma$
such that $\delta(\sigma) = 0$. We assume also that $\Gamma\backslash G/K$
is a spin manifold. Then the range of the trace map is
$$
[\tr] (K_0 ({C}^*_r(\Gamma, \sigma)) ) = \Z\theta + \Z + B, $$
where $2(2\pi)^2\theta = \langle[\omega\cup \omega], [\Gamma]\rangle$ is
the result of
pairing the cup product of
multipliers $[\omega\cup\omega] $ with the fundamental class of $\Gamma$,
and $B= {\rm range}(T_\omega)$. \end{thm}

\noindent{\bf Proof.}

We observe again that by the results of the previous subsection the twisted
Kasparov map is an isomorphism. Therefore to compute the range of the trace
map on $K_0$, it suffices to compute the range of the trace map on elements
of the form $$\mu_\sigma([{\cal E}^0 ] -
[{\cal E}^1])$$
for any element
$$[{\cal E}^0 ] -
[{\cal E}^1] \in K^0(\Gamma\backslash G/K).$$ where ${\cal E}^0, {\cal E}^1$
are vector bundles over the compact manifold $\Gamma\backslash G/K$.

By the twisted $L^2$ index theorem, Theorem 1.7, one has \begin{equation}
[\tr]({\rm ind}_{(\Gamma, \sigma)} (\widetilde{\npartial_{{\cal
E}}^+}\otimes \nabla)) =
\frac{1}{(2\pi)^2}\int_{\Gamma\backslash G/K} \widehat{A}(\Omega)
\tr(e^{R^{{\cal E}}}) e^{\omega}.
\label{(3.6)}
\end{equation}
We next simplify the right hand side of equation (\ref{(3.6)}) using $$
\begin{array}{rl}
\widehat{A}(\Omega) &= 1 -\frac{1}{24} p_1(\Omega)\\ \tr(e^{R^{{\cal E}}})
&= \rank {\cal E} + \tr(R^{{\cal E}}) + \frac{1}{2} \tr({R^{{\cal E}}}^2) \\
e^{\omega} &= 1 + {\omega} + \frac{1}{2} {\omega}^2 . \end{array} $$
Therefore one has
\[
[\tr]({\rm ind}_{(\Gamma, \sigma)} (\widetilde{\npartial_{{\cal
E}}^+}\otimes \nabla)) =
\frac{\rank {\cal E}}{2(2\pi)^2} \; {\displaystyle\int_{\Gamma\backslash
G/K} \omega^2} \]
\[ + \frac{1}{2(2\pi)^2}{ \displaystyle \int_{\Gamma\backslash
G/K}\tr({R^{{\cal E}}}^2)}
-\frac{1}{24(2\pi)^2}{ \displaystyle \int_{\Gamma\backslash G/K}
p_1(\Omega)} +
\frac{1}{(2\pi)^2} { \displaystyle \int_{\Gamma\backslash G/K} \tr(R^{{\cal
E}}) \wedge \omega}, \]
Now by the Atiyah-Singer index theorem \cite{AS}, we see that $$
-\frac{1}{24(2\pi)^2}{ \displaystyle \int_{\Gamma\backslash G/K}
p_1(\Omega)} +
\frac{1}{2(2\pi)^2} { \displaystyle \int_{\Gamma\backslash
G/K}\tr({R^{{\cal E}}}^2)}
= \hbox{index} (\npartial_{{\cal E}}^+)
\in \Z,
$$
Therefore we see that
$$
[\tr] (K_0 ({C}^*_r(\Gamma, \sigma)) ) = \Z\theta + \Z + B, $$
where $2 (2\pi)^2\theta = \displaystyle\int_{\Gamma\backslash G/K}
\omega^2\;$ and
$\quad B = {\rm range}(T_\omega)$.$\diamond$

\begin{rems}{\em
The spin hypothesis on $\Gamma\backslash G/K$ can be easily replaced by
spin$^\C$, without much alteration in the proof.} \end{rems}

The proof of the following proposition is similar to that of Proposition
2.8 and we omit it.

\begin{prop} Let $\Gamma$ be as in Theorem 2.17. Let $\sigma$ be a
multiplier on
$\Gamma$ such that $\delta(\sigma) = 0$. If $\sigma$ defines a rational
cohomology class, then there are at most a finite number of unitary
equivalence classes of projections, other than $0$ and $1$,
in the reduced twisted group $C^*$-algebra ${C}^*_r(\Gamma, \sigma)$.
\label{} \end{prop}

\section{Applications to the spectral theory of projectively periodic
elliptic operators and the classification of twisted group $C^*$ algebras}

In this section, we apply the range of the trace theorem, to prove some
qualitative results on the spectrum of projectively periodic self adjoint
elliptic operators on the universal covering of a good orbifold, or what is
now best known as non-commutative Bloch theory. In particular, we study
generalizations of the hyperbolic analogue of the Ten Martini Problem in
\cite{CHMM} and the Bethe-Sommerfeld conjecture. We also classify up to
isomorphism,
the twisted group $C^*$ algebras for a cocompact Fuchsian group.

Let $M$ be a compact, good orbifold, that is, the universal cover $\Gamma
\to \M \to M$ is a smooth manifold and we will assume as before that there
is a $(\Gamma, \bar\sigma)$-action on $L^2(\M)$ given by $T_\gamma =
U_\gamma \circ S_\gamma \, \forall \gamma \in \Gamma$. Let $\E, \ \F$ be
Hermitian vector bundles on $M$ and let $\E, \ \F$ be the corresponding
lifts to $\Gamma$-invariants Hermitian vector bundles on $\M$. Then there
are $(\Gamma, \sigma)$-actions on $L^2(\M, \E)$ and $L^2(\M, \F)$ which are
also given by $T_\gamma = U_\gamma \circ S_\gamma \, \forall \gamma \in
\Gamma$.

Now let $D : L^2(\M, \E) \to L^2(\M, \F)$ be a self adjoint elliptic
differential
operator that commutes with the $(\Gamma, \bar\sigma)$-action that was
defined earlier.
We begin with some basic facts about the spectrum of such an operator.
Recall that the {\em discrete spectrum} of $D$, $spec_{disc}(D)$ consists
of all the eigenvalues of $D$ that have finite multiplicity, and the {\em
essential spectrum} of $D$, $spec_{ess}(D)$ consists of the complement
$spec(D) \setminus spec_{disc}(D)$. That is, $spec_{ess}(D)$ consists of
the set of accumulation points of the spectrum of $D$, $spec(D)$. Our first
goal is to prove that the essential spectrum is unbounded. Our proof will
be a modification of an argument in \cite{Sar}.

\begin{lemma}
Let $D: L^2(\M, \E) \to L^2(\M, \E)$ be a self adjoint elliptic
differential operator that commutes with the $(\Gamma, \bar\sigma)$-action.
Then the {\em discrete spectrum}
of $D$ is empty.
\label{lemmaS}
\end{lemma}

\noindent{\bf Proof.}
Let $\lambda$ be an eigenvalue of $D$ and $V$ denote the corresponding
eigenspace. Then $V$ is a $(\Gamma, \sigma)$- invariant subspace of
$L^2(\M, \E)$.
If ${\cal F}$ is a relatively compact fundamental domain for the action of
$\Gamma$ on $\M$,
one sees as in section 1.2 that there is a $(\Gamma, \bar\sigma)$-isomorphism
$$
L^2(\M, \E) \cong L^2(\Gamma) \otimes L^2({\cal F}, \E|_{\cal F})$$
Here the
$(\Gamma, \bar\sigma)$-action on $L^2({\cal F}, \E|_{\cal F})$ is trivial,
and is the regular
$(\Gamma, \bar\sigma)$ representation on $L^2(\Gamma)$. Therefore it
suffices to show that
the dimension of {\em any} $(\Gamma, \bar\sigma)$-invariant subspace $V$ of
$L^2(\Gamma)$ is infinite dimensional. Let $\{v_1, \ldots, v_N\}$ be an
orthonormal basis for $V$. Then one has $$ T_\gamma v_i(\gamma')
= \sum_{j=1}^N U_{ij}(\gamma) v_j(\gamma') \qquad \forall \gamma, \gamma'
\in \Gamma
$$
where $U = (U_{ij}(\gamma) )$ is some $N\times N$ unitary matrix. Therefore
$$ \begin{array}{rl}
N= \sum_{j=1}^N ||v_i||^2 & = \sum_{j=1}^N \sum_{\gamma\in \Gamma}
|v_i(\gamma\gamma')|^2\\ &= \sum_{\gamma\in \Gamma} \sum_{i=1}^N
\sum_{j=1}^N \sum_{k=1}^N
U_{ij}(\gamma) \overline{U_{ik}(\gamma)}
v_j(\gamma')\overline{v_k(\gamma')}\\ & = \sum_{\gamma\in \Gamma}
\sum_{j=1}^N |v_j(\gamma')|^2\\ &= \#(\Gamma) \sum_{j=1}^N
|v_j(\gamma')|^2.
\end{array}$$
Since $\#(\Gamma) = \infty$, it follows that either $N=0$ or $N=\infty$.
$\diamond$

\begin{cor}
Let $D$ be as in Lemma \ref{lemmaS} above. Then the {\em essential
spectrum} of $D$ coincides with the spectrum of $D$, and so it is
unbounded. \end{cor}

\noindent{\bf Proof.}
By the Lemma above, we conclude that $spec_{ess}(D)$ and $spec(D)$
coincide. Since $D$ is an unbounded self-adjoint operator, it is a standard
fact that $spec(D)$ is unbounded cf. \cite{Gi}, yielding the result.
$\diamond$

Note that in general the spectral projections of $D$, $E_\lambda \not\in
{C^*} (\sigma)$
(see section 1.2 for the definition). However one has

\begin{prop}[Sunada, Bruning-Sunada]
Let $D$ be defined as in Lemma \ref{lemmaS} above. If $\lambda_0 \not\in
{\hbox{spec}}(D)$, then $E_{\lambda_0} \in {C^*} (\sigma)$. \label{prop6}
\end{prop}

\noindent{\bf Proof.}
Firstly, there is a standard reduction to the case when $D$ is positive and
of even order $d\ge 2$ cf. \cite{BrSu}, so we will assume this without loss
of generality.
By a result of Greiner, see also Bruning-Sunada \cite{BrSu}, there are off
diagonal estimates for the Schwartz kernel of the heat operator $e^{-t D}$
$$
|k_t(x,y)| \le C_1 t^{-n/d}exp\left(-C_2
d(x,y)^{d/(d-1)}t^{-1/(d-1)}\right) $$ for some positive constants $C_1,
C_2$ and for $t>0$ in any compact interval. Since the volume growth of a
orbifold covering space is at most exponential, we see in particular that
$|k_t(x,y)|$ is $L^1$ in both variable separately, so that $$ e^{-t D} \in
{C^*} (\sigma).$$
Note that
$$\chi_{[0,e^{-t\lambda}]}(D) = \chi_{[0,\lambda]}(e^{-t D}).$$ Let
$t=1$ and $\lambda_1 = -\log\lambda_0$. Then $\lambda_1 \not\in
{\hbox{spec}}(e^{- D})$ and
$$
\chi_{[0,\lambda_1]}(e^{- D}) = \phi (e^{- D}) $$ where $\phi$ is a
compactly supported smooth function, $\phi \cong 1$ on $[0,\lambda_1]$ and
$\phi \cong 0$ on the remainder of the spectrum. Since ${C^*} (\sigma)$ is
closed under the continuous functional calculus, it follows that $\phi
(e^{- D}) \in {C^*} ( \sigma)$, that is $E_{\lambda_0} \in {C^*} (\sigma)$.
$\diamond$

Let $D$ be any $(\Gamma, \bar\sigma)$-invariant self-adjoint elliptic
differential operator $D$ on $\M$. Since the spectrum $spec(D)$ is a closed
subset of ${\R}$, its complement ${\R}\backslash spec(D)$ is a countable
union of {\em disjoint} open intervals. Each such interval is called a {\em
gap} in the spectrum.

\begin{prop} Let $\Gamma$ be a Fuchsian group of signature $(g, \nu_1,
\ldots, \nu_n)$.
Let $\sigma$ be a multiplier on
$\Gamma$ such that $\delta(\sigma) = 0$, and $2\pi\theta = \langle[\sigma],
[\Gamma]\rangle \in (0,1]$ be the result of pairing the cohomology class of
$\sigma$ with the fundamental class of $\Gamma$. If $\theta$ is rational,
then
the spectrum of any $(\Gamma, \bar\sigma)$-invariant self-adjoint elliptic
differential operator $D$ on $\M$ has only a finite number of gaps in the
spectrum in every half line $(-\infty, \lambda]$. Here $\Gamma\to \M\to M$
is the universal orbifold covering of a compact good orbifold $M$ with
orbifold fundamental group $\Gamma$.
In particular, the intersection of
$spec(D)$ with any compact interval in ${\R}$ is never a Cantor set.
\label{prop7}
\end{prop}

\noindent{\bf Proof.} We first observe that by equation 1.2 in section 1.2,
one has
${C^*} ( \sigma) \cong {C}^*_r(\Gamma, \sigma) \otimes {\cal K}$.
By Proposition 2.16 and Theorem 2.15, it follows that one has the estimate
$C_\sigma(\Gamma)
\ge 1/q >0$ for the Kadison constant in this case. Then one applies Theorem
1 in Br\"uning-Sunada \cite{BrSu} to deduce the proposition.
$\diamond$

In words, we have shown that
whenever the multiplier is rational, then the spectrum of a projectively
periodic elliptic operator is the union of countably many (possibly
degenerate) closed intervals, which can only accumulate at infinity.

Recall the important $\Gamma$-invariant
elliptic differential operator, which is  the Schr\"odinger operator $$
H_{V} = \Delta +V
$$
where $\Delta$ denotes the Laplacian on functions on $\M$ and $V$ is a
$\Gamma$-invariant function on $\M$. It is known that the Baum-Connes
conjecture is true for all amenable discrete subgroups of a connected Lie
group and also for discrete subgroups of $\SO(n,1)$, see \cite{Kas1} and
$\SU(n,1)$, see \cite{JuKas}. For all these groups $\Gamma$, it follows
that the Kadison constant $C_1(\Gamma)$ is positive. Therefore we see by
the arguments above that the spectrum of the periodic elliptic operator
$H_{V}$ is the union of countably many (possibly degenerate) closed
intervals, which can only accumulate at infinity. This gives evidence for
the following:

\begin{conj}[The Generalized Bethe-Sommerfeld conjecture] The spectrum of
any $\Gamma$-invariant Schr\"odinger operator $H_{V}$ has only a {\em
finite} number of bands, in the sense that the intersection of the
resolvent set with ${\R}$ has only a finite number of components. \end{conj}

We remark that the Bethe-Sommerfeld conjecture has been proved completely
by Skriganov \cite{Skri} in the Euclidean case.

This leaves open the question of whether there are $(\Gamma,
\bar\sigma)$-invariant
elliptic differential operators $D$ on $\J$ with Cantor spectrum when
$\theta$ is irrational. In the Euclidean case, this is usually known as the
{\em Ten Martini Problem}, and is to date, not completely solved, though
much progress has been made \cite{Sh}. We pose a generalization of this
problem to the hyperbolic case (which also includes the Euclidean case):

\begin{conj}[The Generalized Ten Dry Martini Problem] Suppose given a
multiplier $\sigma$ on
$\Gamma$ such that $\delta(\sigma) = 0$, and let $2\pi\theta =
\langle[\sigma], [\Gamma]\rangle \in (0,1]$ be the result of pairing the
cohomology class of $\sigma$ with the fundamental class of $\Gamma$. If
$\theta$ is irrational, then
there is a $(\Gamma, \bar\sigma)$-invariant elliptic differential operator
$D$ on $\J$ which has a Cantor set type spectrum, in the sense that the
intersection of $spec(D)$ with some compact interval in ${\R}$ is a Cantor
set.
\end{conj}

\subsection{The four-dimensional case}

The spectral properties studied in this section do not represent a purely
two-dimensional phenomenon. In fact, it is possible to derive similar
results in higher dimensions, as the following example shows.

\begin{prop}
\label{prop7/8}
Let $\Gamma$ be a discrete, torsion-free cocompact subgroup of $\SO_0(1,
4)$ or of $\SU(1,2)$, and $\sigma$ be a multiplier of $\Gamma$ such that
$\delta(\sigma) = 0$. We assume also that $\Gamma\backslash G/K$ is a spin
manifold. If $[\omega]\in H^2(M,{\R})$ is a rational cohomology class, then
the spectrum of any $(\Gamma, \bar\sigma)$-invariant self-adjoint elliptic
differential operator $D$ on $M$ has only a finite number
of gaps in the spectrum in every half line $(-\infty, \lambda]$. In particular,
the intersection of
$spec(D)$ with any compact interval in ${\R}$ is never a Cantor set. \end{prop}

\noindent{\bf Proof.}
By Proposition 2.19 and Theorem 2.17, it follows that one has the estimate
$C_\sigma(\Gamma)
\ge 1/q >0$ for the Kadison constant in this case. Then one applies Theorem 1
in Br\"uning-Sunada \cite{BrSu} to deduce the proposition.
$\diamond$

\subsection{On the classification of twisted group $C^*$-algebras}
We will now use
the range of the trace theorem \ref{thm5}, to give a complete
classification,
up to isomorphism, of the twisted group $C^*$-algebras ${C}^*(\Gamma,
\sigma)$, where
$\Gamma$ is a Fuchsian group of signature $(g, \nu_1, \ldots, \nu_n)$ and
we assume as before that $\delta(\sigma) = 0$.

\begin{prop}
{\bf (The isomorphism classification of twisted group $C^*$--algebras)} Let
$ \sigma, \sigma' \in H^2(\Gamma, {\R}/ {\Z})$ be multipliers on
$\Gamma$ satisfying $\delta(\sigma) = 0 = \delta(\sigma')$, and $$
2\pi\theta = \langle\sigma, [\Gamma]\rangle \in (0,1], \ \ 2\pi\theta' =
\langle \sigma', [\Gamma]\rangle \in (0,1]$$ be the result of pairing
$\sigma$, $\sigma'$ with the fundamental class of $\Gamma$. Then
${C}^*(\Gamma, \sigma) \cong {C}^*(\Gamma, \sigma')$ if and only if
$$\theta' \in \left\{(\theta + \sum_{i=1}^n \beta_i/\nu_i) \quad
\hbox{mod}\ 1,\ (1-\theta + \sum_{i=1}^n \beta_i/\nu_i) \quad \hbox{mod}\
1\right\},$$ where $0\le \beta_i \le \nu_i-1 \quad \forall i = 1,\ldots
,n$. \label{prop8} \end{prop}

\noindent{\bf Proof.}
Let $\tr$ and $\tr'$ denote the canonical traces on ${C}^*(\Gamma, \sigma)$
and ${C}^*(\Gamma, \sigma')$ respectively. Let $$ \phi : {C}^*(\Gamma,
\sigma) \cong {C}^*(\Gamma, \sigma') $$ be an isomorphism, and let
$$
\phi_* : K_0({C}^*(\Gamma, \sigma)) \cong K_0({C}^*(\Gamma, \sigma')) $$
denote the induced map on $K_0$.
By Theorem \ref{thm5}, the range of the trace map on $K_0$ is $$ [\tr] (K_0
({C}^*(\Gamma, \sigma)) ) = {\Z} \theta + {\Z} + \sum_{i=1}^n
{\Z} (1/\nu_i)
$$
and
$$
[\tr'] (K_0 ({C}^*(\Gamma, \sigma')) ) = {\Z} \theta' +{\Z} +\sum_{i=1}^n
{\Z} (1/\nu_i).
$$
Therefore if $\theta$ is irrational, then $${\Z} \theta + {\Z}
+\sum_{i=1}^n {\Z} (1/\nu_i) = {\Z} \theta' + {\Z} + \sum_{i=1}^n {\Z}
(1/\nu_i) $$ implies that $\theta'$ is also irrational and that $$\theta
\pm \theta' \in {\Z} + \sum_{i=1}^n {\Z} (1/\nu_i) .$$ Since $\theta,
\theta' \in (0,1]$, one deduces that $$\theta' \in \left\{(\theta
+ \sum_{i=1}^n \beta_i/\nu_i)\quad \hbox{mod} 1, (1-\theta + \sum_{i=1}^n
\beta_i/\nu_i) \quad\hbox{mod} 1\right\},$$ where $0\le \beta_i \le \nu_i-1
\quad \forall i = 1,\ldots ,n$. Virtually the same argument holds when
$\theta$ is rational, but one argues in $K$-theory first, and applies the
trace only at the final step.

First observe that a diffeomorphism $C: \Sigma_{g'} \to \Sigma_{g'}$ lifts
to a diffeomorphism $C'$ of ${\bf H}$ such that $C' \Gamma {C'}^{-1} =
\Gamma$, i.e. it defines an automorphism of $\Gamma$. Recall that the
finite group $$
G = \left\{ C_i : \quad C_i^{\nu_i} = 1 \quad \forall i = 1, \ldots, n
\right\} $$
acts on $\Sigma_{g'}$ with quotient $\Sigma(g, \nu_1, \ldots, \nu_n)$. By
the observation above, we see that $G$ also acts as automorphisms of
$\Gamma$. Now $C_i[\Gamma] = \lambda_i[\Gamma]$, where $\lambda_i \in
{\C}$. Since
$C_iC_j[\Gamma] = \lambda_i\lambda_j[\Gamma]$ and $C_i^{\nu_i} = 1$, it
follows that $\lambda_i$ is an $\nu_i^{th}$ root of unity, i.e. $\lambda_i
= e^{2\pi \sqrt{-1} (1/\nu_i)}$. Let $C \in G$, i.e. $C = \prod_{i=1}^n
C_i^{\beta_i}$. We evaluate
$$<C^*[\sigma], [\Gamma]> = <[\sigma], C_*[\Gamma]>
={<[\prod_{i=1}^n\lambda_i^{\beta_i}\sigma], [\Gamma]>} = \theta +
\sum_{i=1}^n \beta_i/\nu_i.$$ As in section 2.1 we see that
$$C^*[\sigma] = [\prod_{i=1}^n\lambda_i^{\beta_i}\sigma]\in \ker\delta
\subset H^2(\Gamma, {\mathbf U}(1)).$$
Therefore the automorphism $C_*$ of $\Gamma$ induces an isomorphism of
twisted group $C^*$-algebras
$$
{C}^*(\Gamma, \sigma) \cong {C}^*(\Gamma, C^*\sigma)\cong {C}^*(\Gamma,
\lambda\sigma).
$$
where $\lambda = \prod_{i=1}^n\lambda_i^{\beta_i}$.

Now let $\psi : \Sigma_{g'} \to \Sigma_{g'}$ be an orientation reversing
diffeomorphism.
Then as observed earlier,
$\psi$ induces an automorphism $\psi_*:\Gamma \to \Gamma$ of $\Gamma$. We
evaluate
$$<\psi^*[\sigma], [\Gamma]> = <[\sigma], \psi_*[\Gamma]>
={\overline{<[\sigma], [\Gamma]>}}=
<[\bar\sigma], [\Gamma]>,$$
since $\psi$ is orientation
reversing. As in section 2.1 we see that $\psi^*[\sigma] =
[\bar\sigma]\in \ker\delta \subset H^2(\Gamma, {\mathbf U}(1))$. Therefore
the automorphism $\psi_*$ of $\Gamma$ induces an isomorphism of twisted
group $C^*$-algebras
$$
{C}^*(\Gamma, \sigma) \cong {C}^*(\Gamma, \psi^*\sigma)\cong {C}^*(\Gamma,
\bar\sigma).
$$

Therefore if
$\theta' \in \{(\theta + \sum_{i=1}^n \beta_i/\nu_i) \quad \hbox{mod} 1,
(1-\theta + \sum_{i=1}^n \beta_i/\nu_i) \quad \hbox{mod} 1\}$, where $0\le
\beta_i \le \nu_i-1 \quad \forall i = 1,\ldots ,n$, one has ${C}^*(\Gamma,
\sigma) \cong {C}^*(\Gamma, \sigma')$, completing the proof of the
proposition.
$\diamond$

\vspace*{12pt}
\noindent{\bf Acknowledgments:} The first author is partially supported by
NSF grant DMS-9802480.  Research by the second author is supported by the
Australian Research Council.

\nonumsection{References}

\end{document}